\documentclass[10pt,twocolumn,twoside]{IEEEtran} 

\usepackage{mathtools}
\usepackage{amsmath,mathrsfs,amsthm}
\usepackage{amssymb}
\usepackage{color}
\usepackage[all]{xy}
\usepackage{graphicx,xspace,bm}
\usepackage{subfigure}
\usepackage{url}
\usepackage[lined,linesnumbered,algoruled]{algorithm2e}
\usepackage{setspace}
\usepackage{url}
\newtheorem{theorem}{Theorem}[section]
\newtheorem{definition}[theorem]{Definition}

\newtheorem{lemma}[theorem]{Lemma}

\newtheorem{remark}[theorem]{Remark}

\newtheorem{corollary}[theorem]{Corollary}
\newtheorem{proposition}[theorem]{Proposition}  

\newcommand{\vertices}{\mathcal{V}}
\newcommand{\edges}{\mathcal{E}}

\newcommand{\Bgraph}{\mathcal{G}}

\newcommand{\BB}{\mathcal{B}}

\newcommand{\JJ}{\mathcal{J}}
\newcommand{\II}{\mathcal{I}}
\newcommand{\KK}{\mathcal{K}}

\newcommand{\GG}{\mathcal{G}}

\newcommand{\HH}{\mathcal{H}}
\newcommand{\NN}{\mathcal{N}}

\newcommand{\EE}{\mathcal{E}}

\newcommand{\zeros}{\mathbf{0}}

\newcommand{\argmax}{\mathrm{argmax}}

\newcommand{\abs}[1]{\ensuremath{\left\lvert{#1}\right\rvert}}
\newcommand{\norm}[1]{\ensuremath{\| #1 \|}}

\newcommand{\real}{{\mathbb{R}}}

\newcommand{\realnonnegative}{{\mathbb{R}}_{\ge 0}}

\newcommand{\integerspositive}{\mathbb{Z}_{\geq 1}}
\newcommand{\integersnonnegative}{\mathbb{Z}_{\geq 0}}

\renewcommand{\min}{\operatorname{min}}
\renewcommand{\max}{\operatorname{max}}
\newcommand{\until}[1]{\{1,\dots,#1\}}
\renewcommand{\until}[1]{[#1]}
\newcommand{\map}[3]{#1:#2 \rightarrow #3}
\newcommand{\setmap}[3]{#1:#2 \rightrightarrows #3}

\newcommand{\setdef}[2]{\{#1 \; | \; #2\}}

\newcommand{\gradient}{\nabla}

\newcommand{\lm}{\lambda}

\newcommand{\baalgo}{\textsc{Bid Adjustment Algorithm}\xspace}
\newcommand{\coeffa}{a}
\newcommand{\coeffc}{c}
\newcommand{\opt}{\mathrm{opt}}

\newcommand{\nout}[1]{\mathcal{N}^{+}_{#1}}
\newcommand{\nin}[1]{\mathcal{N}^{-}_{#1}}
\newcommand{\zov}{\overline{z}}
\newcommand{\yov}{\overline{y}}
\newcommand{\oprocendsymbol}{\hbox{$\bullet$}}
\newcommand{\oprocend}{\relax\ifmmode\else\unskip\hfill\fi\oprocendsymbol}

\newcommand{\longthmtitle}[1]{\mbox{}\textup{\textsl{(#1):}}}

\newcommand{\solone}{\mathrm{Sol}_{\mathrm{sopf}}}
\newcommand{\soltwo}{\mathrm{Sol}_{\mathrm{eff}}}
\newcommand{\nt}{\tilde{n}}
\newcommand{\nb}{\bar{n}}

\renewcommand{\theenumi}{(\roman{enumi})}

\renewcommand{\baselinestretch}{.983}

\newcommand{\myclearpage}{\clearpage}
\renewcommand{\myclearpage}{}

\allowdisplaybreaks

\begin{document}

\title{Iterative bidding in electricity markets:\\ rationality and
  robustness}





\author{Ashish Cherukuri \qquad Jorge Cort\'{e}s\thanks{Preliminary
    versions appeared in conference proceedings
    as~\cite{AC-JC:16-acc} and~\cite{AC-JC:17-ciss}.
  
  Ashish Cherukuri and Jorge~Cort\'{e}s are with the Department of
  Mechanical and Aerospace Engineering, University of California, San
  Diego, \texttt{\{acheruku,cortes\}@ucsd.edu}.}}

\maketitle

\begin{abstract}
  This paper studies an electricity market consisting of an
  independent system operator (ISO) and a group of generators. The
  goal is to solve the DC optimal power flow (DC-OPF) problem: have
  the generators collectively meet the power demand while minimizing
  the aggregate generation cost and respecting line flow limits in the
  network.  The ISO by itself cannot solve the DC-OPF problem as
  generators are strategic and do not share their cost functions.
  Instead, each generator submits to the ISO a bid, consisting of the
  price per unit of electricity at which it is willing to provide
  power.  Based on the bids, the ISO decides how much production to
  allocate to each generator to minimize the total payment while
  meeting the load and satisfying the line limits.  We provide a
  provably correct, decentralized iterative scheme, termed \baalgo, for
  the resulting Bertrand competition game.
  Regarding convergence, we show that the algorithm takes the
  generators’ bids to any desired neighborhood of the efficient Nash
  equilibrium at a linear convergence rate. As a consequence, the
  optimal production of the generators converges to the optimizer of
  the DC-OPF problem. Regarding robustness, we show that the algorithm
  is robust to affine perturbations in the bid adjustment scheme
  and that there is no incentive for any individual generator to
  deviate from the algorithm by using an alternative bid update
  scheme. We also establish the algorithm robustness to collusion,
  i.e., we show that, as long as each bus with generation has a
  generator following the strategy, there is no incentive for any
  group of generators to share information with the intent of tricking
  the system to obtain a higher payoff. Simulations illustrate our results.
\end{abstract}

\section{Introduction}\label{sec:Intro}

As part of the plan to integrate distributed energy resources (DERs)
into the electricity grid, regulating authorities envision a
hierarchical architecture where, at the lower layer, different sets of
DERs coordinate their response under an aggregator and, at the upper
layer, the independent system operator (ISO) interacts with the
aggregators to solve the optimal power flow (OPF) problem.  In this
scenario, aggregators function as (virtual, large-capacity)
generators, and the aggregation would allow DERs to participate into
markets in which, individually, they do not have the capacity to do
so.  While the DERs under one aggregator can cooperate among
themselves, the aggregators compete with each other in the electricity
market.
In this paper, we focus on the competition aspect of this vision: we
study policies that individual generators, in conjunction with the
ISO, can implement to solve the OPF problem while acting in a selfish
and rational fashion.



\subsubsection*{Literature review}
The study of competition in electricity markets is a classical
topic~\cite{SS:02,DK-GS:04}.  Extensively studied models are supply
function, Bertrand (price) and Cournot (capacity) bidding,
see~\cite{RJ-JNT:11,WT-RJ:13-cdc,DC-SB-AW:16-arXiv}, respectively, and
references therein.  These studies analyze the properties of the game
that different bidding models result into by determining the existence
of the Nash equilibrium of the game and estimating its
efficiency. Some
works~\cite{TL-MS:05,XH-DR:07,GIB-CC-MK-FS:10,UN-GP:10}, on the other
hand, propose iterative algorithms for the players that compute the
Nash equilibrium of the game.  However, these algorithms either
require generators to have some information about other generators
(cost functions or bids) or assume that the demand of each generator
is a continuous function of the bids.  Our work does not make any such
assumptions, which also rules out the possibility of using various
other Nash equilibrium learning algorithms, such as
best-response~\cite{EN-RG-RRJ:10}, fictitious play~\cite{DF-DKL:98},
and extremum seeking~\cite{PF-MK-TB:12,MSS-KHJ-DMS:12}.
In a related set of works~\cite{JW-PG-SM-MM:12,DJS-MC-AMA:16},
decentralized generation planning is achieved by assuming the
generators to be price-takers and designing iterative schemes based on
dual-decomposition~\cite{DPP-MC:07}.  In our work, however, we
consider a strategic scenario where generators bid into the market and
are hence price-setters.  The work~\cite{MNF-SD:16} proposes an
iterative auction algorithm for a market where both generators and
consumers are strategic but does not provide convergence guarantees
for the generated bid sequences.  The paper~\cite{SP-AG:16}, closer in
spirit to our work, proposes an iterative method for the generators to
find the Nash equilibrium assuming they do not have any information
about each other.  At each iteration, the generators send to the ISO
the gradient of their cost functions at a certain generation value and
the ISO then adjusts these generation values so that social welfare is
maximized. An important difference between this setup and ours is that
we do not assume truthful bidding of gradient information by the
generators.

Our electricity market game belongs to the broader class of
multi-leader-single-follower games~\cite{SL-TM:10,JSP-MF:05}. The Nash
equilibria of such games can be thought of as optimizers of
mathematical programs with equilibrium constraints
(MPEC)~\cite{ZQL-JSP-DR:08}, 
that are traditionally solved in a centralized
manner~\cite{MCF-SPD-AM:05}.  The work~\cite{DZ-GHL:14} provides a
distributed method to find the equilibria of an MPEC problem but
requires the follower's (the ISO in our case) optimization to have a
unique solution for each action of the leaders (the generators). This
is in general not the case for electricity markets.  Our work broadly
relates to the recent developments in the area of ``learning in
games'', see e.g.,~\cite{SH-AMC:15,TL-VS-EV:16-soda}, and references
therein.  Learning mechanisms proposed in there do not apply directly
to the electricity market setting as they do not consider network
constraints for allocation of goods.  Finally, our work has close
connections with the growing interest in the design of provably
correct distributed algorithms for the cooperative solution of
economic dispatch, see~\cite{AC-JC:16-auto,STC-ADDG-CNH:15-tcst}, and
references therein.

\subsubsection*{Statement of contributions}

The background for the inelastic electricity market game considered
here is that the ISO seeks to find the production levels that solve
the DC optimal power flow (DC-OPF) problem for a group of strategic
generators which do not share their cost functions. Consequently, the
ISO cannot solve the DC-OPF problem by itself.  However, each
generator submits a bid to the ISO specifying the price per unit of
electricity at which the generator is willing to provide power. Given
these bids, the ISO decides how much production to allocate to each
generator so that the cost of generation is minimized, the loads are
met, and the network flow constraints are satisfied. The resulting
Bertrand competition model defines the game among the generators,
where the actions are the bids and the payoffs are the profits.  We
define the concept of the efficient Nash equilibrium, that is, the Nash
equilibrium at which the generators are willing to produce the amount
that corresponds to the optimizer of the DC-OPF problem.
Our first contribution gives two set of conditions that ensure
existence and uniqueness, respectively, of an efficient Nash
equilibrium for the inelastic electricity market game.  Our second
contribution is the design of the \baalgo along with its correctness
analysis.  This algorithm can be understood as ``learning via repeated
play'', where generators are ``myopically selfish'', changing their
bid at each iteration with the aim of maximizing their own payoff.
Along the execution, the only information available to the generators
is their bid and the amount of generation that the ISO request from
them.  In particular, generators are not aware of the number of other
generators, their costs, bids, or payoffs.  We show that this
decentralized iterative scheme is guaranteed to take the bids of the
generators to any neighborhood of the unique efficient Nash
equilibrium provided the stepsizes are chosen appropriately. Further,
we establish that the convergence rate is linear. Our third
contribution analyzes the robustness properties of the
\baalgo. Specifically, we establish that the convergence is not affected
by affine disturbances, thus showing that deviations in stepsizes by
the generators can be handled gracefully. Additionally, we show that
there is no incentive for any individual generator to deviate from the
algorithm by using an alternative bid update scheme.  Finally, we also
show that, if at each generator bus there is at least one generator
running the \baalgo, then there is no incentive for other generators
connected to the network to not follow the algorithm, i.e., this
adjustment scheme becomes a rational choice for all generators. These
properties provide a sound justification for why the group generators
would adopt this iterative bid adjustment scheme to solve the DC-OPF
problem. Simulations illustrate our results.
%


%
%
\subsubsection*{Notation}
Let $\real$, $\realnonnegative$, $\integersnonnegative$,
$\integerspositive$ be the set of real, nonnegative real, nonnegative
integer, and positive integer numbers, respectively. We use the shorthand
notation $\until{n}$ to denote the set $\{1,\dots,n\}$.  The $2$-norm on
$\real^n$ is represented by $\norm{\cdot}$.  Let $\BB_{\delta}(x) =
\setdef{y \in \real^n}{\norm{y-x} < \delta}$ be the open ball centered
at $x \in \real^n$ with radius $\delta > 0$. 
Given $x,y\in \real^n$, $x_i$ is the $i$-th component
of $x$, and $x \le y$ denotes $x_i \le y_i$ for $i \in \until{n}$. We
use $\zeros_N = (0,\ldots,0) \in \real^N$.
We let
$[u]^{+} = \max \{0,u\}$ for~$u \in \real$.
A \emph{directed graph}, or simply
\emph{digraph}, is represented by a pair $\Bgraph=(\vertices,\edges)$, where
$\vertices$ is the vertex set and $ \edges \subseteq \vertices\times
\vertices $ is the edge set.  
For a digraph, $\nout{v_i} = \setdef{v_j \in \vertices}{(v_i, v_j) \in
  \edges}$ and $\nin{v_i} = \setdef{v_j \in \vertices}{(v_j,v_i) \in
  \edges}$ are the sets of out- and in-neighbors of $v_i$,
respectively.

\section{Problem statement}\label{sec:problem}

Consider an electrical power network with $N_b \in \integerspositive$
buses.  The physical interconnection between the buses is given
by a digraph $\GG = (\vertices, \edges)$, where nodes
correspond to buses and edges to physical power lines.  The direction
for each edge represents the convention of positive power flow.  The
power flow on the line $(i,j) \in \EE$ is $z_{ij} \in \real$.  Each
line $(i,j) \in \EE$ has a limit on the power flowing through it
(in either direction), represented by $\zov_{ij} > 0$.  Assume that
each bus $i \in \until{N_b}$ is connected to $n_i \in
\integersnonnegative$ strategic generators.
We let $N = \sum_{i=1}^{N_b} n_i$ be the total number of generators
and assign them a unique identity in $\until{N}$. Let the set of
generators at node $i$ be $G_i \subset \until{N}$ (this set is empty
if there are no generators connected to bus $i$).  The power demand at
bus $i$ is denoted by $y_i \ge 0$ and is assumed to be fixed and known
to the Independent System Operator (ISO) that acts as the central
regulating authority.  The total demand is $\yov =
\sum_{i=1}^{N_b} y_i$.
The cost $f_n(x_n)$ of generating $x_n \in
\realnonnegative$ amount of power by the $n$-th generator is given by
a quadratic function
\begin{align}\label{eq:cost-quadratic}
  f_n(x) = \coeffa_n x^2 + \coeffc_n x,
\end{align}
where $\coeffa_n > 0$ and $\coeffc_n \ge 0$. 
Given a power allocation $x = (x_1, \dots, x_{N}) \in
\realnonnegative^{N}$, the aggregate cost
is $\sum_{n=1}^{N} f_n(x_n)$.  The \emph{dc optimal power flow
  problem} (DC-OPF) consists of
\begin{subequations}\label{eq:edp-flow}
  \begin{align}
    \underset{(x,z)}{\text{minimize}} & \quad \sum_{n=1}^{N} f_n(x_n),
    \label{eq:edp-flow-0}
    \\
    \text{subject to} & \quad \sum_{j \in \NN_i^+} \!\! z_{ij} -
    \sum_{j \in \NN_i^-} \!\!z_{ij} = \sum_{n \in G_i} \! x_n - y_i,
    \; \forall i, \label{eq:edp-flow-1}
    \\
    & \quad -\zov_{ij} \le z_{ij} \le \zov_{ij}, \, \, \forall
    (i,j), \label{eq:edp-flow-2}
    \\
    & \quad x \ge \zeros_N. \label{eq:edp-flow-3}
  \end{align}
\end{subequations}
This problem finds the generation profile that meets the load at each
bus (ensured by~\eqref{eq:edp-flow-1}), respects the line
constraints (due to~\eqref{eq:edp-flow-2}), and minimizes the total
cost (given by the objective
function~\eqref{eq:edp-flow-0}). In~\eqref{eq:edp-flow-1} we make the
convention that if $G_i = \emptyset$, then the first term on the
right-hand side is zero.  We assume that~\eqref{eq:edp-flow} is
feasible. Since the individual costs are quadratic, the
optimizer of the problem, denoted~$(x^*,z^*)$, is
unique~\cite{SB-LV:04}.

The goal for the ISO is to solve~\eqref{eq:edp-flow}.  The ISO can
interact with the generators, whereas each generator can only
communicate with the ISO and is not aware of the number of other
generators participating in the market and their respective cost
functions, or the load at its own bus.  While the ISO knows the loads
and the limits on the power lines, it does not have any information
about the cost functions of the generators.  Therefore, power
allocation is decided following a bidding process, resulting into a
game-theoretic formulation. Instead of sharing their cost with the
ISO, the generators bid the price per unit of power that they are
willing to provide the power at. This price-based bidding is
well known in the economics literature as Bertrand
competition~\cite[Chapter 12]{AMC-MW-JG:95}.  Specifically, generator
$n$ bids the cost per unit power $b_n \in \realnonnegative$ and, when
convenient, we denote the bids of all other generators except $n$
by $b_{-n} = (b_1, \dots, b_{n-1}, b_{n+1}, \dots, b_N)$.  Given 
the bids $b = (b_1, \dots, b_N) \in \realnonnegative^N$, the ISO
solves the following \emph{strategic dc optimal power flow problem}
(S-DC-OPF)
\begin{subequations}\label{eq:sedp-flow}
  \begin{align}
    \underset{(x,z)}{\text{minimize}} & \quad  \sum_{n=1}^N b_n x_n,
    \label{eq:sedp-flow-0}
    \\
    \text{subject to} & \quad  \sum_{j \in \NN_i^+} z_{ij} -
    \sum_{j \in \NN_i^-} z_{ij} =  \sum_{n \in G_i} x_n - y_i, \, \,
    \forall i, \label{eq:sedp-flow-1}
    \\
    & \quad -\zov_{ij} \le z_{ij} \le \zov_{ij}, \, \, \forall (i,j),
    \label{eq:sedp-flow-2}
    \\
    & \quad x \ge \zeros_N. \label{eq:sedp-flow-3}
  \end{align}
\end{subequations}
The difference between~\eqref{eq:sedp-flow}
and~\eqref{eq:edp-flow} is the objective function which is linear in the
former and nonlinear, convex in the latter.  The ISO
solves~\eqref{eq:sedp-flow} once all the bids are gathered. Let
$(x^{\opt}(b),z^{\opt}(b))$ be the optimizer of~\eqref{eq:sedp-flow}
that the ISO selects (note that there might not be a unique optimizer)
given bids $b$. This determines the power requested from each generator,
given by the vector $x^{\opt}(b)$. Knowing this process, the objective
of each generator $n$ is to bid a quantity $b_n \ge 0$ that maximizes
its payoff $\map{u_n}{\realnonnegative^2}{\real}$,
\begin{align}\label{eq:payoff}
  u_n(b_n,x_n^{\opt}(b)) = b_n x_n^{\opt}(b) - f_n(x_n^{\opt}(b)),
\end{align}
where $x_n^{\opt}(b)$ is the $n$-th component of the 
optimizer $x^{\opt}(b)$. 

\begin{definition}\longthmtitle{Inelastic electricity market
    game}\label{def:inelastic-game}
  The \emph{inelastic electricity market} game is defined by the
  following
    \begin{enumerate}
    \item Players: the set of generators $\until{N}$,
    \item Action: for each player $n$, the bid $b_n \in
      \realnonnegative$,
    \item Payoff: for each player $n$, the payoff $u_n$
      in~\eqref{eq:payoff}.
    \end{enumerate}
\end{definition}

Wherever convenient, for any $n \in \until{N}$, we use interchangeably
the notation $b$ and $(b_n,b_{-n})$, as well as, $x^{\opt}(b)$ and
$x^{\opt}(b_n,b_{-n})$.  Note that the payoff of the players is not
only defined by the bids of other players but also by the optimizer
of~\eqref{eq:sedp-flow} that the ISO selects. For this reason, the
definition of the pure Nash equilibrium for the game described below
is slightly different from the standard one, see e.g.~\cite{DF-JT:91}.

\begin{definition}\longthmtitle{Nash equilibrium}\label{def:Nash}
  The \emph{(pure) Nash equilibrium} of the inelastic electricity
  market game is the bid profile of the group $b^* \in
  \realnonnegative^N$ for which there exists an optimizer
  $(x^{\opt}(b^*),z^{\opt}(b^*))$ of the
  optimization~\eqref{eq:sedp-flow} that satisfies
  \begin{equation}\label{eq:nash}
    u_n(b_n,x_n^{\opt}(b_n,b_{-n}^*)) \le u_n(b_n^*,x_n^{\opt}(b^*)) ,
  \end{equation}
  for all $n \in \until{N}$, all bids $b_n \in \realnonnegative$, and
  all optimizers $(x^{\opt}(b_n,b_{-n}^*),z^{\opt}(b_n,b_{-n}))$
  of~\eqref{eq:sedp-flow} given bids $(b_n,b_{-n}^*)$.
\end{definition}

We are specifically interested in bid profiles for which the optimizer
of the DC-OPF problem is also a solution to the S-DC-OPF problem. This
is captured in the following definition.

\begin{definition}\longthmtitle{Efficient bid}\label{def:efficient}
  An \emph{efficient bid} of the inelastic electricity market is a bid
  $b^* \in \realnonnegative^N$ for which the optimizer $(x^*,z^*)$
  of~\eqref{eq:edp-flow} is also an optimizer of~\eqref{eq:sedp-flow}
  given bids $b^*$ and
  \begin{equation}\label{eq:efficiency}
    x_n^* = \argmax_{x \ge 0} b_n^* x - f_n(x), \quad \text{ for all } n
    \in \until{N}.
  \end{equation}
\end{definition}

The right-hand side of~\eqref{eq:efficiency} is unique as costs are quadratic.

\begin{definition}\longthmtitle{Efficient Nash
    equilibrium}\label{def:efficient-Nash}
  A bid $b^*$ is an \emph{efficient Nash equilibrium} of the inelastic
  electricity market game if it is an efficient bid and is a Nash
  equilibrium.
\end{definition}

At the efficient Nash equilibrium, the production that the generators
are willing to provide, maximizing their profit, coincides with the
optimal generation for the DC-OPF problem~\eqref{eq:edp-flow}. This
property justifies the study of efficient Nash equilibria. Note that
given the efficient bid profile, there might be many solutions
to~\eqref{eq:sedp-flow} because the problem is linear. Thus, 
the ISO might not be able to find $x^*$
given the efficient bid.  However, once the ISO knows that an
efficient Nash equilibrium bid is submitted, it can ask the generators
to also submit the desirable generation levels at that bid, which
would exactly correspond to the solution of the DC-OPF problem.

\section{Existence and uniqueness of efficient Nash
equilibrium}\label{sec:existence}

Here, we establish the existence of an efficient Nash
equilibrium of the inelastic electricity market game described in
Section~\ref{sec:problem} and provide a condition for its uniqueness.

\begin{proposition}\longthmtitle{Existence of efficient Nash
    equilibrium}\label{pr:existence}
  Assume that at each bus of the network either there is more than one
  generator or there is none, i.e., either $n_i = 0$ or $n_i \ge 2$ for
  each $i \in \until{N_b}$.  Then, there exists an efficient Nash
  equilibrium of the inelastic electricity market game.
\end{proposition}
\begin{IEEEproof}
  For convenience, we write~\eqref{eq:edp-flow-1}
  and~\eqref{eq:edp-flow-2} as
  \begin{align*}
    J_1 z - J_2 x + y = 0 \quad \text{ and } \quad J_3 z \le
    \overline{z}_c,
  \end{align*}
  respectively. Here, $J_1 \in \{0,1,-1\}^{N_b \times N_b}$ defines 
  the interconnection of buses in the digraph $\GG$,
  specifically, $(i,j)$-th
  element of $J_1$ is $1$ if $(i,j) \in \EE$, is $-1$ if $(j,i) \in
  \EE$, and $0$ otherwise. The matrix $J_2 \in \{0,1\}^{N_b \times N}$
  defines the connectivity of generators to buses, that is, $(i,j)$-th
  element of $J_2$ is $1$ if and only if $j$-th generator is connected
  to $i$-th bus. Lastly, 
  \begin{align*}
    J_3 = \begin{bmatrix} I_{\abs{\EE}} \\ -I_{\abs{\EE}} \end{bmatrix}
    \quad \text{ and } \quad \overline{z}_c = \begin{bmatrix}
      \overline{z} \\ \overline{z}
  \end{bmatrix}.
  \end{align*}
  The Lagrangian of the optimization~\eqref{eq:edp-flow} is
  \begin{align*}
    L(x,z,\nu,\mu,\lm) & = \textstyle \sum_{n=1}^N f_n(x_n) + \nu^\top
    ( J_1 z - J_2 x + y)
    \\
    & \qquad + \mu^\top(J_3 z - \overline{z}_c) - \lm^\top x,
  \end{align*}
  where $\nu \in \real^{N_b}$, $\mu \in \realnonnegative^{2
    \abs{\EE}}$, and $\lm \in \realnonnegative^N$ are Lagrange
  multipliers corresponding to
  constraints~\eqref{eq:edp-flow-1},~\eqref{eq:edp-flow-2},
  and~\eqref{eq:edp-flow-3}, respectively. Since
  constraints~\eqref{eq:edp-flow-1}-\eqref{eq:edp-flow-2} are affine
  and the feasibility set is nonempty, the refined Slater condition is
  satisfied for~\eqref{eq:edp-flow} and hence, the duality gap between
  the primal and the dual optimization problems is
  zero~\cite{SB-LV:04}.  Under this condition, a primal-dual optimizer
  $(x^*,z^*,\nu^*,\mu^*,\lm^*)$ satisfies the following
  Karush-Kuhn-Tucker (KKT) conditions
  \begin{subequations}\label{eq:KKT}
    \begin{align}
      \gradient f(x^*) - J_2^\top \nu^* - \lm^* = 0, \label{eq:KKT-1}
      \\
      J_1^\top \nu^* - J_3^\top \mu^* = 0, \label{eq:KKT-2}
      \\
      J_1 z^* - J_2 x^* + y = 0, 
      \\
      J_3 z^* \le \overline{z}_c, \quad x^* \ge 0,
      \\
      \lm^* \ge 0, \quad \mu^* \ge 0, \label{eq:KKT-5}
      \\
      (x^*)^\top \lm^* = 0, \text{ and } (\mu^*)^\top(J_3 z^* -
      \overline{z}_c) = 0, \label{eq:KKT-6}
    \end{align}
  \end{subequations}
  where $\gradient f(x^*) = (\gradient f_1(x^*_1), \gradient
  f_2(x^*_2), \dots, \gradient f_N(x^*_N))^\top$. In the rest of the
  proof, we show that the following bid profile, constructed from a
  primal-dual optimizer, is an efficient Nash equilibrium of the
  inelastic electricity market game
  \begin{align}\label{eq:nash-eq}
    b^*_n =
    \begin{cases}
      \nu^*_{i(n)}, & \quad \text{ if } \min\setdef{x_m^*}{m \in
	G_{i(n)}}
      > 0,
      \\
      \gradient f_n(0), & \quad \text{ otherwise,}
    \end{cases}
  \end{align}
  where $i(n) \in \until{N_b}$ denotes the bus of the network to which
  generator $n$ is connected to.  Given the
  form~\eqref{eq:cost-quadratic} of the cost functions, we deduce $b^*
  \ge 0$. Moreover, from the definition of $J_2$, one can deduce that
  either all generators $n \in G_i$ have $b^*_n = \nu_i$ or all of
  them have $x^*_n = 0$.  Next, to show that the bid
  $b^*$ defined in~\eqref{eq:nash-eq} is efficient, we first establish
  \begin{equation}\label{eq:bid-eff}
    x_n^* = \argmax_{x \ge 0} b_n^* x - f_n(x),
  \end{equation}
  for all $n \in \until{N}$. For each $n$, consider the optimization
  $\max_{x \ge 0} b_n^* x - f_n(x)$.  Because zero duality holds for
  this optimization, a point $x_o \in \realnonnegative$ is an
  optimizer if and only if it satisfies the KKT conditions
  \begin{align*}
    b_n^* -\gradient f_n(x_o) - \mu_o = 0,
    \\
    \mu_o \ge 0, \quad x_o \ge 0, \quad \mu_o x_o = 0,
  \end{align*}
  where $\mu_o$ is the dual optimizer. Since $x_n^*$ satisfies the
  above conditions with $\mu_o = \lm^*_n$, the
  expression~\eqref{eq:bid-eff} holds. To claim the efficiency of
  $b^*$, we next show that $(x^*,z^*)$ is one of the optimizers
  of~\eqref{eq:sedp-flow} given bids $b^*$. Note that the KKT
  conditions for~\eqref{eq:sedp-flow} are given by~\eqref{eq:KKT} with
  the term $\gradient f(x^*)$ in~\eqref{eq:KKT-1} replaced with
  $b^*$. Also, one can show using the KKT conditions~\eqref{eq:KKT}
  and the definition of $b^*$ that $b^* - J_2^\top \nu^* \ge 0$. Using
  these facts, we deduce that $(x^*,z^*,\nu^*,\mu^*,b^* - J_2^\top
  \nu^*)$ satisfies the KKT conditions for~\eqref{eq:sedp-flow} and
  hence, $(x^*,z^*)$ is an optimizer of~\eqref{eq:sedp-flow}.

  Our final step is to show the Nash equilibrium
  condition~\eqref{eq:nash} for the bid profile $b^*$.  Note that for
  each $n$, the payoff at the bid profile-optimizer pair
  ($b^*,x^{\opt}(b^*)) = (b^*,x^*)$ is nonnegative. Specifically, if
  $x^*_n = 0$, then $u_n(b_n^*,x_n^{\opt}(b^*)) = 0$.  If $x^*_n > 0$,
  using the fact that $\gradient f_n(x) \le b_n^*$ for all $x \in
  [0,x_n^*]$, we get
  \begin{align*}
    u_n(b_n^*,x_n^{\opt}(b^*)) & = b_n^* x_n^* - f_n(x_n^*)
    \\
    & = \int_0^{x_n^*} \gradient (b_n^* x - f_n(x)) dx
    \\
    & = \int_0^{x_n^*} (b_n^* - \gradient f_n(x)) dx \ge 0.
   \end{align*}
   Now pick any generator $n \in \until{N}$.  For bid $b_n \not =
   b_n^*$ we have two cases, first, $b_n > b_n^*$ and second, $b_n \le
   b_n^*$.  For the first case, either (i) $x^*_n = 0$ which implies
   that $x^{\opt}(b_n,b_{-n}^*) = 0$ and so $u_n(b_n,
   x_n^{\opt}(b_n,b_{-n}^*)) = u_n(b_n^*,x_n^*) = 0$; or (ii) $x^*_n >
   0$, so all bids at bus $i(n)$ are equal, implying that $n$
   increasing its bid yields $x_n^{\opt}(b_n,b_{-n}^*) = 0$. That is,
   $u_n(b_n, x_n^{\opt}(b_n,b_{-n}^*)) = 0 \le u_n(b_n^*,x_n^*)$.  For
   the second case,
   \begin{align*}
     u_n(b_n,x_n^{\opt}(b_n,b_{-n}^*)) & = b_n
     x_n^{\opt}(b_n,b_{-n}^*) - f_n(x_n^{\opt}(b_n,b_{-n}^*))
     \\
     & \le b_n^* x_n^{\opt}(b_n,b_{-n}^*) -
     f_n(x_n^{\opt}(b_n,b_{-n}^*))
     \\
     & \le b_n^* x_n^* - f_n(x_n^*) = u_n(b_n^*,x_n^*),
   \end{align*}
   where in the first inequality we use $b_n \le b_n^*$ and in
   the second we use~\eqref{eq:bid-eff}. This
   shows~\eqref{eq:nash}, concluding the
   proof.
\end{IEEEproof}

Note that the condition in Proposition~\ref{pr:existence} of having
zero or at least two generators at each bus is reasonable. If this is
not the case, i.e., there is a bus with a single generator, and the
line capacities are such that the load at that bus can only be met by
that generator, then there is possibility of market manipulation.  The
generator at the bus can set its bid arbitrarily high as no other
generator can meet that load and consequently, there does not exist a
Nash equilibrium.  Next we provide a sufficient condition that ensures
uniqueness of the efficient bid. 

\begin{lemma}\longthmtitle{Uniqueness of the efficient
    bid}\label{le:uniqueness-1}
  Assume that the optimizer $x^*$ of~\eqref{eq:edp-flow} satisfies
  $x^*_n > 0$ for all $n \in \until{N}$. Then, there exists a unique
  efficient bid $b^* \in \realnonnegative^N$ of the inelastic
  electricity market game given by
  \begin{align}\label{eq:eff-bid}
    b^*_n = \gradient f_n(x^*_n) = 2 \coeffa_n x^*_n + \coeffc_n, \quad
    \text{ for all } n
  \end{align}
\end{lemma}
\begin{IEEEproof}
  By definition, an efficient bid $b \in \realnonnegative^N$ satisfies 
  \begin{align*}
    x^*_n = \argmax_{x \ge 0} b_n x - f_n(x)
  \end{align*}
  for all $n$. Since $x^*_n > 0$, first-order optimality condition of
  the above optimization yields $b_n = \gradient f_n(x^*_n)$. This
  establishes~\eqref{eq:eff-bid} and hence, the uniqueness. 
\end{IEEEproof}

From Proposition~\ref{pr:existence} and Lemma~\ref{le:uniqueness-1},
we conclude the following result.

\begin{corollary}\longthmtitle{Uniqueness of the efficient Nash
    equilibrium}\label{cor:uniqueness}
  Assume that at each bus of the network either there is more than one
  generator or there is none. Further assume that 
  the optimizer $x^*$ of~\eqref{eq:edp-flow} satisfies $x^*_n
  > 0$ for all $n \in \until{N}$. Then, there exists a unique
  efficient Nash equilibrium of the inelastic electricity market game
  given by~\eqref{eq:eff-bid} for all $n$. 
\end{corollary}

In the rest of the paper, we assume that the sufficient conditions in
Corollary~\ref{cor:uniqueness} hold unless otherwise stated. Note that
the definition of the unique efficient Nash equilibrium given
in~\eqref{eq:eff-bid} is consistent with the one provided
in~\eqref{eq:nash-eq}. This is so because if $x^*_n > 0$ for all $n$,
then $\gradient f_{\tilde{n}} (x^*_{\tilde{n}}) = \nu^*_i$
for each bus $i \in \until{N_b}$ and every generator $\tilde{n} \in
G_i$.

\section{The \baalgo and its convergence properties}\label{sec:algo}

In this section, we introduce a decentralized Nash equilibrium seeking
algorithm, termed \baalgo. We show that its executions lead the
generators to the unique efficient Nash equilibrium, and consequently,
to the optimizer of the DC-OPF problem~\eqref{eq:edp-flow}.

\subsection{\baalgo}

We start with an informal description of the \baalgo.  The algorithm
is iterative and can be interpreted as ``learning via repeated play''
of the inelastic electricity market game by the generators.  Both ISO
and generators have bounded rationality, with each generator trying to
maximize its own profit and the ISO trying to maximize the welfare of
the entities.
\begin{quote}
  \emph{[Informal description]:} At each iteration $k$, generators
  decide on a bid and send it to the ISO. Once the ISO has obtained
  the bids, it computes an optimizer of the S-DC-OPF
  problem~\eqref{eq:sedp-flow}, denoted $(x^{\opt}(k),z^{\opt}(k))$
  for convenience, and sends the corresponding production level at the
  optimizer to each generator. At the $(k+1)$-th iteration, generators
  adjust their bid based on their previous bid, the amount of produced
  power that maximizes their payoff for the previous bid, and the
  allocation of generation assigned by the ISO. The iterative process
  starts with the generators arbitrarily selecting initial bids that
  yield a positive profit.
\end{quote}

The \baalgo is formally presented in
Algorithm~\ref{alg:inelastic-two}.


\begin{algorithm}
  \SetAlgoLined
  \DontPrintSemicolon
  \SetKwFor{Case}{case}{}{endcase}
  \SetKwInOut{giv}{Data} \SetKwInOut{ini}{Initialize}
  \SetKwInOut{state}{State} \SetKwInput{start}{Initiate}
  \SetKwInOut{msg}{Messages} \SetKwInput{Kw}{Executed by}
  \Kw{generators $n \in \until{N}$ and ISO}
  \giv{cost $f_n$ and stepsizes $\{\beta_k\}_{k \in \integerspositive}$ for each
    generator $n$, and load $y$ for ISO}
  \ini{Each generator $n$ selects arbitrarily $b_n(1) \ge \coeffc_n$, sets
    $k=1$, and jumps to step~\ref{step1}; ISO sets $k=1$ and
    waits for step~\ref{step2}}
  \BlankLine \While{$k>0$}{ \tcc{For each generator $n$:} Receive
  $x^{\opt}_n(k-1)$ from ISO \; 
  Set $b_n(k) \! = \! [b_n(k \! - \! 1) \!+\! \beta_k(x^{\opt}_n(k \! -
    \! 1) - q_n(k \! - \! 1))]^+$ \; \label{step3}
    Set $q_n(k) = \argmax_{q \ge 0} b_n(k) q -
    f_n(q)$ \; \label{step4}
    \label{step1} Send $b_n(k)$ to the ISO; set $k = k+1$ \;
    \BlankLine \tcc{For ISO:} Receive $b_n(k)$ from each $n  \in
      \until{N}$
    \;
    \label{step2} Find a solution $(x^{\opt}(k),z^{\opt}(k))$
    to~\eqref{eq:sedp-flow} given $b(k)$ \;
    %
    Send $x^{\opt}_n(k)$ to each  $n \in \until{N}$; set $k = k+1$}
  \caption{\baalgo}\label{alg:inelastic-two}
\end{algorithm}
In the \baalgo, the role of the ISO is to compute an optimizer of the
S-DC-OPF problem after the bids are submitted.  Generators adjust
their bids at each iteration in a ``myopically selfish'' and rational
fashion, with the sole aim of maximizing their payoff.
Intuitively,
\begin{LaTeXdescription}
  \item[if $n$ gets $x_n^{\opt}(k) = 0$,] two things can happen: (i) $n$ was
  willing to produce a positive quantity $q_n(k) > 0$ at bid $b_n(k)$
  but the demand from the ISO is $x^{\opt}_n(k) = 0$. Thus, the rational choice
  for $n$ is to decrease the bid in the next iteration and
  increase its chances of getting a positive payoff; (ii)
  $n$ was willing to produce nothing $q_n(k) = 0$ at $b_n(k)$ and
  got $x^{\opt}_n(k)=0$. At this point, reducing the bid will not increase
  the payoff as it will not be willing to produce more at a lower
  bid. Alternatively, increasing the bid will not make the amount
  that the ISO wants the generator to produce positive. 
  Hence, the bid stays put.
\item[if $n$ gets $x_n^{\opt}(k) > 0$,] then it would want to move the bid
  in the direction that makes its payoff higher in the next iteration,
  assuming that $n$ gets a positive generation signal from the ISO in
  the next round. If $q_n(k) <
  x_n^{\opt}(k)$, then the demand from the ISO is more than what the
  generator is willing to produce, so $n$ increases its cost, i.e.,
  the bid. If $q_n(k) > x_n^{\opt}(k)$, then the demand is less than what
  the generator is willing to supply so $n$ decreases its bid.
\end{LaTeXdescription}

\begin{remark}\longthmtitle{Information structure and other
    learning approaches}\label{re:algo-comparison}
  {\rm Generators have no knowledge of the number of other players,
    their actions, or their payoffs.  The only information they have 
    at each iteration is their own bid and the generation that
    the ISO requests from them. This information structure rules out
    the applicability of a number of Nash equilibrium learning
    methods, including best-response dynamics~\cite{EN-RG-RRJ:10},
    fictitious play~\cite{DF-DKL:98}, or other gradient-based
    adjustments~\cite{GIB-CC-MK-FS:10}, all requiring some kind of
    information about other players. Methods that relax this
    requirement, such as extremum seeking used
    in~\cite{PF-MK-TB:12,MSS-KHJ-DMS:12}, rely on the payoff functions
    being continuous in the actions of the players, which is not the
    case for the inelastic electricity market game.  \oprocend }
\end{remark}

\begin{remark}\longthmtitle{Stopping criteria and justification of
    ``myopically selfish'' strategies}\label{re:myopic}
  {\rm Algorithm~\ref{alg:inelastic-two} consists of 
  an infinite number of iterations. To make it implementable, 
  later we identify stopping criteria, see
  Remark~\ref{re:stop-criteria}, based on a parameter that the ISO
  selects. Since this is not known to the generators, they cannot
  predict when the algorithm will terminate and, hence, they do not have
  an incentive to play strategically to maximize their payoff in the
  long term. Given this, they should focus on maximizing the payoff in
the next iteration, which justifies the myopically selfish perspective
adopted here. \oprocend }
\end{remark}

%

\subsection{Convergence analysis}
In this section, we show that the generator bids along any execution
of the \baalgo converge to a neighborhood of the unique efficient Nash
equilibrium.  The size of the neighborhood is a decreasing function of
the stepsize and can be made arbitrarily small.

We first present a series of results that highlight certain geometric
properties of the bid update done in Step~\ref{step3} of
Algorithm~\ref{alg:inelastic-two}. These results form the basis for
establishing later the convergence guarantee.  The following result
states that one could neglect the projection operator in
Step~\ref{step3} of Algorithm~\ref{alg:inelastic-two}.

\begin{lemma}\longthmtitle{Generator bids are lower
    bounded}\label{le:seq-prop}
  In Algorithm~\ref{alg:inelastic-two}, let $0 < \beta_k < 2
  \coeffa_n$ for all $n \in \until{N}$ and $k \in
  \integerspositive$. Then, $b_n(k) \ge \coeffc_n$ and for all $n \in
  \until{N}$ and $k \in \integerspositive$, 
  \begin{align}\label{eq:q-func}
    q_n(k) = \frac{b_n(k) - \coeffc_n}{2 \coeffa_n}.
  \end{align}
\end{lemma}
\begin{IEEEproof}
  Equation~\eqref{eq:q-func} follows directly from $b_n(k) \ge
  \coeffc_n$, so we focus on proving the latter.  We proceed by
  induction. Note that $b_n(1) \ge \coeffc_n$ for all $n \in
  \until{N}$. Assume that $b_n(k) \ge \coeffc_n$ for some $k \in
  \integerspositive$ and let us show $b_n(k+1) \ge \coeffc_n$.  We
  have
  \begin{align*}
    \! b_n(k+1) \! &= [b_n(k) + \beta_k(x_n^{\opt}(k) - q_n(k))]^+
    \\
    & \overset{(a)}{\ge} [b_n(k) \! - \! \beta_k q_n(k)]^+ \!
    \overset{(b)}{=} \!  \Bigl[ b_n(k) \! - \! \beta_k \Bigl( \frac{b_n(k) -
      \coeffc_n}{2 \coeffa_n} \Bigr) \Bigr]^+
    \\
    & = \Bigl[ \Bigl(1 - \frac{\beta_k}{2 \coeffa_n} \Bigr) b_n(k) + \beta_k
    \frac{\coeffc_n}{2 \coeffa_n} \Bigr]^+
    \\
    & \overset{(c)}{=} \Bigl(1 - \frac{\beta_k}{2 \coeffa_n} \Bigr) b_n(k) +
    \beta_k \frac{\coeffc_n}{2 \coeffa_n},
  \end{align*}
  where $(a)$ is due to the fact that $x^{\opt}_n(k) \ge 0$, $(b)$ follows
  from the definition of $q_n(k)$ given the fact that $b_n(k) \ge
  \coeffc_n$, and $(c)$ follows from the assumption that $\beta_k < 2
  \coeffa_n$ for all $n$ (which makes both terms in the expression
  positive).  By contradiction, assume $b_n(k+1) < \coeffc_n$. Then,
  \begin{align*}
    \Bigl( 1 - \frac{\beta_k}{2 \coeffa_n} \Bigr) b_n(k) + \beta_k
    \frac{\coeffc_n}{2 \coeffa_n} < \coeffc_n,
  \end{align*}
  which implies that $b_n(k) < \coeffc_n$, a contradiction.
\end{IEEEproof}

Our next result gives a different expression for the bid update step
(cf. Step~\ref{step3}) presenting a geometric perspective of the
direction along which the bids are moving. Specifically, we write the
$k+1$-th bid as the addition of two vectors. The first one is a convex
combination of the $k$-th bid and the efficient Nash equilibrium $b^*$.
Hence, the first vector is closer to $b^*$ as
compared to the $k$-th bid. The second one depends on the difference
between what the ISO requests from the generators and the optimizer
of~\eqref{eq:edp-flow}. If the second term is small enough, then we are
assured that the bids move towards $b^*$. 


\begin{lemma}\longthmtitle{Geometric characterization of the bid
    update}\label{lem:b-int}
  In Algorithm~\ref{alg:inelastic-two}, let $0 < \beta_k < 2
  \coeffa_n$ for all $n \in \until{N}$ and $k \in
  \integerspositive$. Then, we have
  \begin{align*}
    b(k+1) = b^{\operatorname{coc}}(k+1) + \beta_k (x^{\opt}(k) - x^*),
  \end{align*}
  for all $k \in \integerspositive$, where for each $n \in \until{N}$,
  \begin{align*}
    b_n^{\operatorname{coc}}(k+1) = \Bigl(1 - \frac{\beta_k}{2 \coeffa_n}
    \Bigr) b_n(k) + \frac{\beta_k}{2 \coeffa_n} b_n^* .
  \end{align*}
\end{lemma}
\begin{IEEEproof}
  In the proof of Lemma~\ref{le:seq-prop}, we have shown that for all
  $n$ and $k$, the term
  inside the projection operator $[\cdot]^+$ in Step~\ref{step3} of
  Algorithm~\ref{alg:inelastic-two} is nonnegative. Hence, the
  projection can be dropped and we can write
  \begin{align*}
    b_n(k+1) & = b_n(k) + \beta_k ( x_n^{\opt}(k) - q_n(k))
    \\
    & \overset{(a)}{=} b_n(k) + \beta_k (x_n^{\opt}(k)) - \beta_k \Bigl(
    \frac{b_n(k) - \coeffc_n}{2 \coeffa_n} \Bigr)
    \\
    & = \Bigl(1 - \frac{\beta_k}{2 \coeffa_n} \Bigr) b_n(k) + \beta_k
    \Bigl(x_n^{\opt}(k) + \frac{\coeffc_n}{2 \coeffa_n} \Bigr)
    \\
    & = \Bigl(1 - \frac{\beta_k}{2 \coeffa_n} \Bigr) b_n(k) + \beta_k
    (x_n^{\opt}(k) - x_n^*)
    \\
    & \qquad \qquad + \beta_k \Bigl(x_n^* + \frac{\coeffc_n}{2
      \coeffa_n} \Bigr)
    \\
    & \overset{(b)}{=} \Bigl(1 - \frac{\beta_k}{2 \coeffa_n} \Bigr)
    b_n(k) + \frac{\beta_k}{2 \coeffa_n} b_n^* + \beta_k (x_n^{\opt}(k) -
    x_n^*).
  \end{align*}
  In the above expression, we have used~\eqref{eq:q-func} in the
  equality $(a)$ and~\eqref{eq:eff-bid} 
  in the equality $(b)$.
\end{IEEEproof}

The next result gives a lower bound on the inner product between the
direction in which the bids move and the direction towards the
efficient Nash equilibrium.

\begin{lemma}\longthmtitle{Bids move in the direction of the efficient
    Nash equilibrium}\label{le:bound-angle}
  In Algorithm~\ref{alg:inelastic-two}, let $0 < \beta_k < 2
  \coeffa_n$ for all $n \in \until{N}$ and $k \in
  \integerspositive$. Let $\coeffa_{\max} = \max_{n}
  \{\coeffa_n\}$. Then, for all $k \in \integerspositive$,
  \begin{align}\label{eq:term1-bound}
    \langle b(k+1) - b(k), b^* - b(k) \rangle \ge  \frac{\beta_k}{2
      \coeffa_{\max}} \norm{b(k) - b^*}^2.
  \end{align}
\end{lemma}
\begin{IEEEproof}
  Using Lemma~\ref{lem:b-int}, we write
  \begin{align*}
    &\langle b(k+1) - b(k), b^* - b(k) \rangle
    \\
    & \quad = \langle b(k+1) - b^{\operatorname{coc}}(k+1), b^* - b(k)
    \rangle
    \\
    & \qquad \qquad \qquad \qquad + \langle
    b^{\operatorname{coc}}(k+1) - b(k), b^* - b(k) \rangle
    \\
    & \quad = \beta_k \langle x^{\opt}(k) - x^*, b^* - b(k) \rangle +
    \sum_{n=1}^N \frac{\beta_k}{2 \coeffa_n} (b_n^* - b_n(k))^2
    \\
    & \quad \textstyle \overset{(a)}{\ge} \sum_{n=1}^N \frac{\beta_k}{2
      \coeffa_n} (b_n^* - b_n(k))^2 \ge \frac{\beta_k}{2
      \coeffa_{\max}} \norm{b(k) - b^*}^2.
  \end{align*}
  For the inequality $(a)$, we have used the fact that
  \begin{align*}
    \langle x^{\opt}(k) - x^*, b^* - & b(k) \rangle = \bigl(\langle
    x^{\opt}(k),
    b^* \rangle - \langle x^* , b^* \rangle \bigr)
    \\
    & + \bigl(\langle x^* , b(k) \rangle - \langle x^{\opt}(k), b(k)
    \rangle \bigr) \ge 0.
  \end{align*}
  The last inequality follows from the fact that $x^*$ and
  $x^{\opt}(k)$ are the optimizers of~\eqref{eq:sedp-flow} given $b^*$
  and $b(k)$, resp., making both expressions on the right-hand
  side nonnegative.
\end{IEEEproof}

The next result states that the distance between consecutive bids
decreases as the bids get closer to $b^*$. 
In combination with Lemma~\ref{le:bound-angle}, one can see
intuitively that the bids get closer to $b^*$ and, as they get closer
to it, the bid update step behaves as if the bids are reaching an
equilibrium of the update scheme. These two facts lead to convergence.

\begin{lemma}\longthmtitle{Distance between consecutive bids is
    upper bounded}\label{le:bound-dist}
  In Algorithm~\ref{alg:inelastic-two}, let $0 < \beta_k < 2 a_n$ for
  all $n \in \until{N}$ and $k \in \integerspositive$.  Let
  $\coeffa_{\min} = \min_{n} \{\coeffa_n\}$.  Then, for all $k \in
  \integerspositive$,
  \begin{align}\label{eq:term2-bound}
    \norm{b(k+1) - b(k)}^2 \le \frac{\beta_k^2}{2 \coeffa_{\min}^2}
    \norm{b(k) - b^*}^2 + 8 \beta_k^2 \yov^2.
  \end{align}
\end{lemma}
\begin{IEEEproof}
  Consider the following 
  \begin{align}
    & \norm{b(k+1) - b(k)}^2 \notag
    \\
    & \quad \overset{(a)}{=} \sum_{n=1}^N \Bigl( \frac{\beta_k}{2
    \coeffa_n} (b_n^* - b_n(k)) + \beta_k (x_n^{\opt}(k) - x_n^*)
    \Bigr)^2 \notag
    \\
    & \quad \overset{(b)}{\le} \sum_{n=1}^N 2 \Bigl(\frac{\beta_k}{2
      \coeffa_n} (b^*_n - b_n(k)) \Bigr)^2 + \sum_{n=1}^N 2 \beta_k^2
      (x^{\opt}_n(k) - x^*_n)^2 \notag
    \\
    & \quad \overset{(c)}{\le} \frac{\beta_k^2}{2 \coeffa_{\min}^2}
    \norm{b(k) - b^*}^2 + 2 \beta_k^2 \norm{x^{\opt}(k) - x^*}^2.
    \label{eq:diff-ineq}
  \end{align}
  In the above expression, (a) follows from the expression of
  $b_n(k+1)$ from Lemma~\ref{lem:b-int}, (b) follows from the
  inequality $(x+y)^2 \le 2 (x^2 + y^2)$ for $x,y \in \real$, and (c)
  follows from the definition of $\coeffa_{\min}$. Note that
  \begin{align*}
    \norm{x^{\opt}(k) - x^*} \le \sum_{n=1}^N \abs{x^{\opt}_n(k) - x^*_n} \le
    \sum_{n=1}^N \abs{x^{\opt}_n(k)} + \abs{x^*_n} 
    \\
    \textstyle = \sum_{n=1}^N (x^{\opt}_n(k) +
    x^*_n) = 2 \yov.
  \end{align*}
  The proof concludes by using the above bound
  in~\eqref{eq:diff-ineq}.
\end{IEEEproof}

We are ready to present the main convergence result.

\begin{theorem}\longthmtitle{Convergence of the
    \baalgo}\label{th:convergence}
  In Algorithm~\ref{alg:inelastic-two}, let $0 < \beta_k < 2
  \coeffa_n$ for all $n \in \until{N}$ and $k \in
  \integerspositive$. Further, let $0 < r < \norm{b(1) - b^*}$ and for
  all $k \in \integerspositive$ assume
  \begin{align}\label{eq:beta-cond-3}
    \alpha \le \beta_k \le B(r) := \frac{1}{2 a_{\max}} \Bigl(\frac{1}{2
      a_{\min}^2} + \frac{16 \yov^2}{r^2} \Bigr)^{-1},
  \end{align}
  for some $\alpha > 0$. Then, the following holds
  \begin{enumerate}
  \item there exists $l \in \integerspositive$ such that $\norm{b(l) -
      b^*} < r$ and for all $k \in \until{l-1}$, we have $\norm{b(k) -
      b^*} \ge r$ with
    \begin{align}\label{eq:rate}
      \norm{b(k+1) - b^*} \le \Bigl(1 - \frac{\alpha}{2
        \coeffa_{\max}} \Bigr)^{k/2} \norm{b(1) - b^*},
    \end{align}
  \item for all $k \ge l$, 
    \begin{align}\label{eq:ultimate-bound}
      \norm{b(k) - b^*} \le \Bigl(1 + \frac{B(r)}{2 \coeffa_{\max}}
      \Bigr)^{1/2} r.
    \end{align}
  \end{enumerate}
\end{theorem}
\begin{IEEEproof}
  Assume that $\norm{b(k) - b^*} \ge r$ for some $k \in
  \integerspositive$. Then, the upper bound on the stepsizes in 
  the inequality~\eqref{eq:beta-cond-3}
  holds when $r$ is replaced with $\norm{b(k) - b^*}$, that is, $\beta_k
  \le B(\norm{b(k)-b^*})$ for all $k \in \integerspositive$. This is
  because $r \mapsto B(r)$ is strictly increasing in the domain $r > 0$.
  Proceeding with this replacement and
  reordering~\eqref{eq:beta-cond-3}, we obtain
  \begin{align*}
    \beta_k \Bigl( \frac{\norm{b(k) - b^*}^2}{2 a_{\text{min}}^2} + 16
    \yov^2 \Bigr) \le \frac{1}{2a_{\text{max}}} \norm{b(k) - b^*}^2,
  \end{align*}
  or equivalently,  
  \begin{align}\label{eq:beta-cond-4}
    \frac{\beta_k}{2 a_{\min}^2} \norm{b(k) - b^*}^2 + 16 \beta_k
    \yov^2 - \frac{1}{a_{\max}} \norm{b(k) - b^*}^2 \notag
    \\
    \le - \frac{1}{2 a_{\max}} \norm{b(k) - b^*}^2.
  \end{align}
  Now consider the following inequalities
  \begin{subequations}\label{eq:iter-bound-all}
    \begin{align}
      \norm{b(k+1) & - b^*}^2 = \norm{b(k+1) - b(k) + b(k) - b^*}^2
      \notag
      \\
      & = \norm{b(k+1) - b(k)}^2 + \norm{b(k) - b^*}^2 \notag
      \\
      & \qquad + 2 \langle b(k+1) - b(k), b(k) - b^* \rangle \notag
      \\
      & \overset{(a)}{\le} \frac{\beta_k^2}{2 a_{\min}^2} \norm{b(k) -
        b^*}^2 + 8 \beta_k^2 \yov^2 + \norm{b(k) - b^*}^2 \notag
      \\
      & \qquad \qquad - \frac{\beta_k}{ a_{\max}} \norm{b(k) - b^*}^2
      \label{eq:iter-bound-0}
      \\
      & \overset{(b)}{\le} \Bigl(1- \frac{\beta_k}{2 a_{\max}} \Bigr)
      \norm{b(k) - b^*}^2, \label{eq:iter-bound}
    \end{align}
  \end{subequations}
  %
  %
  where in (a) we have used the bounds~\eqref{eq:term1-bound}
  and~\eqref{eq:term2-bound} from Lemmas~\ref{le:bound-angle}
  and~\ref{le:bound-dist}, respectively, and the inequality (b) is
  implied by that in~\eqref{eq:beta-cond-4}. Note that the
  inequality~\eqref{eq:beta-cond-4} is conservative in the sense that
  the term $16\beta_k \yov^2$ could be replaced with $8\beta_k \yov^2$
  and the inequality~\eqref{eq:iter-bound} would still follow. However,
  we opt for this conservativeness while defining the map $r \mapsto
  B(r)$ in~\eqref{eq:beta-cond-3} because it results into robustness
  guarantees for the algorithm as discussed in the forthcoming
  section. Therefore,~\eqref{eq:iter-bound}
  holds whenever $\norm{b(k) - b^*} \ge r$. By assumption, we have $0
  < \Bigl( 1 - \frac{\beta_k}{2 \coeffa_{\max}} \Bigr) < 1$,
  $\norm{b(1) - b^*} > r$, and $\beta_k \ge \alpha$ for all $k \in
  \integerspositive$. Using these facts and
  applying~\eqref{eq:iter-bound} recursively, we conclude part~(i).

  For part (ii), note that if $\norm{b(k) - b^*} \ge r$ for some $k
  \ge l$, then $\norm{b(k+1) - b^*} < \norm{b(k) - b^*}$
  by~\eqref{eq:iter-bound}. Therefore, to find an upper bound on
  $\norm{b(k) - b^*}$ for all $k \ge l$, we only need to consider the
  case when $\norm{b(k) - b^*} < r$. Plugging this bound
  in~\eqref{eq:iter-bound-0} and neglecting the negative term, we get
  \begin{align}\label{eq:norm-r-bound}
    \norm{b(k+1) - b^*}^2 \le \frac{\beta_k^2 r^2}{2 \coeffa_{\min}^2}
    + 8 \beta_k^2 \yov^2 + r^2.
  \end{align}
  From~\eqref{eq:beta-cond-3}, we have
  \begin{align*}
    \frac{\beta_k^2 r^2}{2 \coeffa_{\min}^2} +
    16 \beta_k^2 \yov^2 \le \frac{\beta_k r^2}{2 \coeffa_{\max}}.
  \end{align*}
  The result now follows by upper bounding the right-hand side
  of~\eqref{eq:norm-r-bound} with the left-hand side of the above
  expression and then employing the bound on the stepsizes give
  in~\eqref{eq:beta-cond-3}.
  %
  %
\end{IEEEproof}

\begin{remark}\longthmtitle{Convergence properties from
    Theorem~\ref{th:convergence}}\label{re:conv-prop}
    {\rm The assertion (i) of Theorem~\ref{th:convergence} implies that
    for any choice of $r>0$, one can select stepsizes according
    to~\eqref{eq:beta-cond-3} so that bids reach the set $\BB_r(b^*)$
    in a finite number of steps and at a linear rate.
    Further, once bids reach the set $\BB_{r}(b^*)$, we are assured
    from assertion (ii) that they remain in a neighborhood of $b^*$,
    where the size of the neighborhood is proportional to $r$
    (cf.~\eqref{eq:ultimate-bound}).  In combination, the above facts
    mean that bids converge to any neighborhood of the efficient Nash
    equilibrium at a linear rate provided the stepsizes are selected
    appropriately.  Note that as $r$ becomes small, $B(r)$ gets small
    and so does~$\alpha$. Thus, from~\eqref{eq:rate}, the rate of
    convergence decreases as $r$ becomes small. This presents a
    trade-off between the desired precision and the rate of
    convergence.  \oprocend }
\end{remark}

\begin{remark}\longthmtitle{Stopping criteria for the
    ISO}\label{re:stop-criteria}
  {\rm From the proof of Theorem~\ref{th:convergence}(i) note that, as
    long as $\norm{b(k) - b^*} > r$, the distance to the efficient
    Nash equilibrium decreases. Therefore, if $\norm{b(k) - b^*} > r$
    and $k < l$, then one can write
    \begin{align}
      \norm{b(k+1) & - b(k)} = \norm{b(k+1) - b^* + b^* - b(k)} \notag
      \\
      & \ge \norm{b(k) - b^*} - \norm{b(k+1) - b^*} \notag
      \\
      & \overset{(a)}{\ge} \norm{b(k) - b^*} - \Bigl(1-
      \frac{\alpha}{2 \coeffa_{\max}} \Bigr)^{1/2} \norm{b(k) - b^*}
      \notag
      \\
      & = \Bigl(1 - \Bigl(1- \frac{\alpha}{2
        \coeffa_{\max}}\Bigr)^{1/2} \Bigr) \norm{b(k) - b^*},
	\label{eq:dist-b}
    \end{align} 
    where in (a) we have used~\eqref{eq:iter-bound} and
    $\beta_k \ge \alpha$.  Given this observation, if the ISO has an
    estimate of $\alpha$ and $\coeffa_{\max}$, then it can design a
    stopping criteria based on the distance between consecutive
    bids. In fact, if the ISO decides selects $\epsilon > 0$ and stops
    the iteration whenever $\norm{b(k+1) - b(k)} \le \epsilon$, then
    it has the guarantee that either of the following is
    satisfied
    \begin{enumerate}
      \item the condition $\norm{b(k) - b^*} > r$ and $k < l$ is met and 
	from~\eqref{eq:dist-b} we get 
	\begin{align}\label{eq:stop-bound}
          \norm{b(k) - b^*} \le \epsilon \Bigl(1 - \Bigl(1-
          \frac{\alpha}{2 \coeffa_{\max}}\Bigr)^{1/2} \Bigr)^{-1};
        \end{align}
      \item $\norm{b(k) - b^*} \le r$; or 
      \item $k > l$ in which case from~\eqref{eq:ultimate-bound} we get
	\begin{align*}
	  \norm{b(k) - b^*} \le \Bigl(1 + \frac{B(r)}{2 \coeffa_{\max}}
	  \Bigr)^{1/2} r.
        \end{align*}
    \end{enumerate}
    The ISO does not know the value of $r$; its value depends on the
    stepsizes that the generators select. Assuming that stepsizes are
    small, 
    the ISO can adjust $\epsilon$ depending on the desired accuracy
    level to get the guarantee~\eqref{eq:stop-bound} for the $k$-th
    bid.  For small $\epsilon$, the stopping criteria might never be
    met if stepsizes are too~big.  \oprocend }
\end{remark}


\section{Robustness of the \baalgo}\label{sec:robust}

Here we study the robustness properties of the \baalgo in a
variety of scenarios. We first show that the introduction of
disturbances in the bid update mechanism does not destroy the
algorithm convergence properties.  We then study robustness
against either an individual agent or colluding agents changing their
strategy to get a higher payoff.

\subsection{Robustness to disturbances}\label{sec:robust-dist}

Here we establish the robustness properties of the \baalgo in the
presence of disturbances by characterizing its input-to-state
stability (ISS) properties~\cite{ZPJ-YW:01}.  Let
$\map{d}{\integerspositive}{\real^N}$ model the disturbance to the bid
update mechanism. Such disturbances might arise from agents using
different stepsizes than the prescribed one or other disruption to the
prescribed bid update scheme.  The resulting perturbed version of the
\baalgo can be written as the following discrete-time dynamical system
\begin{subequations}\label{eq:ddyn}
  \begin{align}
    b(k+1) & = [b(k) + \beta_k (x^{\opt}(k) - q(k))+ d(k)]^+,
    \label{eq:ddyn-1}
    \\
    x^{\opt}(k+1) & \in \solone(b(k+1)), \label{eq:ddyn-2}
    \\
    q(k+1) & = \soltwo(b(k+1)), \label{eq:ddyn-3}
  \end{align}
\end{subequations}
where $\setmap{\solone}{\realnonnegative^N}{\realnonnegative^N}$ and
$\map{\soltwo}{\realnonnegative^N}{\realnonnegative^N}$ map a bid
profile to the set of optimizers of problem~\eqref{eq:sedp-flow}
and~\eqref{eq:efficiency}, respectively. Note that $\solone$ is a
set-valued map since~\eqref{eq:sedp-flow} is a linear program.  If $d
\equiv 0$, then the dynamics~\eqref{eq:ddyn} represents the $k$-th
iteration of the \baalgo.

The next result shows that the perturbed version of the
algorithm~\eqref{eq:ddyn} retains the convergence properties of the
unperturbed version provided the magnitude of the disturbance
satisfies an upper bound dependent on the state of the bid.
 
\begin{proposition}\longthmtitle{The \baalgo is robust to
    perturbations in the
    bid update}\label{pr:iss}
  For dynamics~\eqref{eq:ddyn}, let the hypotheses of
  Theorem~\ref{th:convergence} hold and assume that $b_n(k) \ge
  \coeffc_n$ for all $n \in \until{N}$ and $k \in \integerspositive$.
  Let $0 < \theta < \frac{1}{6} \bigl( 1 - \frac{\alpha}{2
    \coeffa_{\max}} \bigr)$ and assume $\norm{d(k)} \le \theta
  \norm{b(k) - b^*}$ for all $k \in \integerspositive$.  Then, the
  following holds
  \begin{enumerate}
  \item there exists $l \in \integerspositive$ such that $\norm{b(l) -
      b^*} < r$ and, for all $k \in \until{l-1}$, we have $\norm{b(k) -
      b^*} \ge r$ with
    \begin{align}\label{eq:rate-perturbed}
      & \norm{b(k+1) - b^*} \notag
      \\
      & \quad \le \Bigl(1 - \frac{\alpha}{2 \coeffa_{\max}} +2 \theta
      + 4 \theta^2 \Bigr)^{k/2} \norm{b(1) - b^*},
    \end{align}
  \item for all $k \ge l$, 
    \begin{align}\label{eq:ultimate-bound-perturbed}
      \norm{b(k) - b^*} \le \Bigl(1 + \frac{B(r)}{2 \coeffa_{\max}} +
      2 \theta + 4 \theta^2 \Bigr)^{1/2} r.
    \end{align}
  \end{enumerate}
\end{proposition}
\begin{IEEEproof}
  Since $b_n(k) \ge \coeffc_n$, we obtain for dynamics~\eqref{eq:ddyn},
    $q_n(k) = \frac{b_n(k) - \coeffc_n}{2 \coeffa_n}$,
  for all $n \in \until{N}$ and $k \in \integerspositive$. Moreover,
  mimicking Lemma~\ref{lem:b-int}, we 
  rewrite the bid update~\eqref{eq:ddyn-1} as 
  \begin{equation}\label{eq:ddyn-bint-update}
    b(k+1) =  b^{\operatorname{coc}}(k+1) + \beta_k (x^{\opt}(k) -
    x^*) + d(k),
  \end{equation}
  for all $k \in \integerspositive$.  Using~\eqref{eq:ddyn-bint-update}
  and following the steps of Lemma~\ref{le:bound-angle} for
  dynamics~\eqref{eq:ddyn-1} we get,
  \begin{align}
    \langle b(k+1) - b(k), & b(k) - b^* \rangle  \le \langle d(k), b(k)
    - b^* \rangle  \notag
    \\
    & \qquad \qquad - \frac{\beta_k}{2 \coeffa_{\max}} \norm{b(k) - b^*}^2,
    \label{eq:ddyn-angle-bound}
  \end{align}
  for all $n \in \until{N}$ and $k \in \integerspositive$. Similarly,
  from the reasoning of Lemma~\ref{le:bound-dist} we obtain 
  \begin{align}
	&\norm{b(k+1) - b(k)}^2 \notag
	\\
	& \le \frac{\beta_k^2}{2 \coeffa_{\min}^2}
	\norm{b(k) - b^*}^2 + 2 \Bigl(\norm{\beta_k(x^{\opt}(k)  - x^*) +
	d(k))} \Bigr)^2  \notag
	\\
	& \le \frac{\beta_k^2}{2 \coeffa_{\min}^2}
	\norm{b(k) - b^*}^2 + 4 \beta_k^2 \norm{x^{\opt}(k) - x^*}^2 + 4
	\norm{d(k)}^2 \notag
	\\
	& \le  \frac{\beta_k^2}{2 \coeffa_{\min}^2}
	\norm{b(k) - b^*}^2 + 16 \beta_k^2 \yov^2 + 4 \norm{d(k)}^2
	\label{eq:ddyn-dist-bound}
  \end{align}
  for all $k \in \integerspositive$ for dynamics~\eqref{eq:ddyn-1}.
  Employing~\eqref{eq:ddyn-angle-bound} and~\eqref{eq:ddyn-dist-bound},
  assuming $\norm{b(k) - b^*} \ge r$, 
  and writing the set of inequalities~\eqref{eq:iter-bound-all} with
  $\alpha \le \beta_k$,
  we deduce the following
  \begin{align}
    \norm{b(k+1) - b^*}^2 & \le \Bigl(1- \frac{\alpha}{2 \coeffa_{\max}} \Bigr) 
    \norm{b(k) - b^*}^2 + 4 \norm{d(k)}^2 \notag
    \\
    & \qquad + 2 \langle d(k) , b(k) - b^* \rangle.
    \label{eq:greater-r-d-ineq}
  \end{align}
  Finally, using $\norm{d(k)} \le \theta \norm{b(k) - b^*}$ we get
  \begin{equation}\label{eq:greater-r-theta-ineq}
    \norm{b(k+1) - b^*}^2  \le \Bigl( 1- \frac{\alpha}{2 \coeffa_{\max}} + 2
    \theta + 4 \theta^2 \Bigr) \norm{b(k) - b^*}^2.
  \end{equation}
  Iteratively, we obtain~\eqref{eq:rate-perturbed}. The
  bound~\eqref{eq:ultimate-bound-perturbed} can be computed in a
  similar way as done in the proof of Theorem~\ref{th:convergence}. 
\end{IEEEproof}

Similar to the convergence guarantees of Theorem~\ref{th:convergence},
the above result establishes that the perturbed version of the
algorithm~\eqref{eq:ddyn} converges to a neighborhood of the efficient
Nash equilibrium provided the stepsizes and the disturbance satisfy
appropriate bounds, and that the size of this neighborhood is tunable
as a function of these.

%
%
%
The next result complements Proposition~\ref{pr:iss} by giving an
alternative representation of robustness of~\eqref{eq:ddyn}. It
establishes two properties: first, when the disturbance is bounded
(not necessarily satisfying the bound of Proposition~\ref{pr:iss}),
the bids remain bounded; second, if the disturbance goes to zero, then
the bids satisfy the bound~\eqref{eq:ultimate-bound} asymptotically.
Notice that both these results do not follow directly from
Proposition~\ref{pr:iss}, justifying the need for a formal proof.

\begin{proposition}\longthmtitle{Bounded disturbance implies bounded
    bids for \baalgo}\label{pr:bounded}
  For dynamics~\eqref{eq:ddyn}, let the hypotheses of
  Theorem~\ref{th:convergence} hold and assume that $b_n(k) \ge
  \coeffc_n$ for all $n \in \until{N}$ and $k \in \integerspositive$.
  Let $\norm{d(k)} \le d_{\max}$ for all $k \in \integerspositive$ and
  let $\theta \in \Bigl(0, \frac{1}{6} \Bigl( 1 - \frac{\alpha}{2
    a_{\max}} \Bigr) \Bigr)$. Then, the following holds for all $k \in
    \integerspositive$,
  \begin{align}
    \norm{b(k) - b^*} & \le \Bigl(1 - \frac{\alpha}{2 a_{\max}} + 2 \theta
    + 4 \theta^2 \Bigr)^{k/2} \norm{b(1) - b^*} \notag
    \\
    & \qquad \qquad \qquad + G(r,\theta,d_{\max}), \label{eq:iss-bound}
  \end{align}
  where $G(r,\theta,d_{\max}) := \max\{ G_1(r,d_{\max}),
  G_2(\theta,d_{\max}) \}$ and 
  \begin{align*}
    G_1(r,d_{\max}) & := \Bigl( \frac{B(r) r^2}{2 a_{\max}} + (2
    d_{\max} + r)^2 \Bigr)^{1/2},
    \\
    G_2(\theta,d_{\max}) & := \Bigl(2 +
    \frac{1}{\theta} \Bigr) d_{\max}.
  \end{align*}
  As a consequence, as $k \to \infty$, if $\norm{d(k)} \to
  0$, then 
  \begin{align}\label{eq:d-zero-limit}
    \max \Bigl\{ \norm{b(k)-b^*}, \Bigl(1 + \frac{B(r)}{2 \coeffa_{\max}}
    \Bigr)^{1/2} r \Bigr\} \to 0.
  \end{align}
\end{proposition}
\begin{IEEEproof}
  We first show that if for some $k \in \integerspositive$,
  $\norm{b(k) - b^*} \le G(r,\theta,d_{\max})$, then $\norm{b(l) -
    b^*} \le G(r, \theta, d_{\max})$ for all $l \ge k$. To this end,
  as a first case, assume that $r \le \norm{b(k) - b^*} \le
  G(r,\theta,d_{\max})$. Then, following the steps of the proof of
  Proposition~\ref{pr:iss}, we arrive
  at~\eqref{eq:greater-r-d-ineq}. If $\norm{d(k)} \le \theta
  \norm{b(k) - b^*}$, then we get the
  inequality~\eqref{eq:greater-r-theta-ineq} which implies that
  $\norm{b(k+1) - b^*} \le \norm{b(k) - b^*} \le
  G(r,\theta,d_{\max})$. On the other hand, if $\norm{d(k)} > \theta
  \norm{b(k) - b^*}$, then using this bound
  in~\eqref{eq:greater-r-d-ineq}, we get
  \begin{align*}
    \norm{b(k+1) - b^*}^2 & < \Bigl(1 - \frac{\alpha}{2 a_{\max}} \Bigr)
    \frac{\norm{d(k)}^2}{\theta^2} + 4 \norm{d(k)}^2 
    \\
    & \qquad \qquad + 2 \frac{\norm{d(k)}^2}{\theta}
    \\
    & < \Bigl( \frac{1}{\theta^2} + \frac{4}{\theta} + 4 \Bigr)
    \norm{d(k)}^2.
  \end{align*}
  Thus, using $\norm{d(k)} \le d_{\max}$, we get $\norm{b(k+1) - b^*}
  < G_2(\theta,d_{\max}) \le G(r,\theta,d_{\max})$. As a second case,
  assume $\norm{b(k) - b^*} < r$. Note that $r <
  G(r,\theta,d_{\max})$, and so $\norm{b(k) -b^*} <
  G(r,\theta,d_{\max})$. For this case, using $\norm{b(k) - b^*} < r$
  and 
  inequalities~\eqref{eq:ddyn-angle-bound}
  and~\eqref{eq:ddyn-dist-bound}, we get as 
  in~\eqref{eq:iter-bound-0} that
  \begin{align*}
    \norm{b(k+1) - b^*}^2 & \le \frac{\beta_k^2 r^2}{2 a_{\min}^2} + 16
    \beta_k^2 \yov^2 + 4 \norm{d(k)}^2 + r^2 
    \\
    & \qquad \qquad + 2 r \norm{d(k)}.
  \end{align*}
  Now applying bounds $\norm{d(k)} \le d_{\max}$ and $\beta_k \le B(r)$,
  we obtain $\norm{b(k+1) - b^*} \le G_1(r,d_{\max})$. Hence, we arrive
  at the conclusion that if
  $\norm{b(k) - b^*} \le G(r,\theta,d_{\max})$, then $\norm{b(l) - b^*}
  \le G(r,\theta,d_{\max})$ for all $l \ge k$. 

  Consider now the case when for some $k \in \integerspositive$,
  $\norm{b(k) - b^*} > G(r,\theta,d_{\max})$. By definition of
  $G(r,\theta,d_{\max})$, this implies that $\norm{b(k) - b^*} > r$ and
  $\norm{b(k) - b^*} > \frac{d(k)}{\theta}$. Therefore, from the proof
  of Proposition~\ref{pr:iss}, we arrive
  at~\eqref{eq:greater-r-theta-ineq}. Finally, combining the reasoning
  of the two cases when $\norm{b(k) - b^*}$ is greater than or less than
  equal to $G(r,\theta,d_{\max})$, we obtain the
  inequality~\eqref{eq:iss-bound}. The limit~\eqref{eq:d-zero-limit}
  follows from that fact that as $k \to \infty$, the first term
  of~\eqref{eq:iss-bound} converges to zero and as $d_{\max}$ tends to
  zero, $G(r,\theta,d_{\max})$ tends to $\Bigl(1 + \frac{B(r)}{2
    a_{\max}} \Bigr)^{1/2} r$.
\end{IEEEproof}

One can observe from~\eqref{eq:iss-bound} that the limiting behavior
of the bids depend on the magnitude of $r$ and $d_{\max}$: if $r$ is
designed to be small enough and if $d_{\max}$ is small enough, or this
bound becomes small as the algorithm iterates, then the bids do
converge to a small neighborhood of $b^*$. 

As an aside, in the theory of ISS for discrete-time dynamical
systems~\cite{ZPJ-YW:01}, one typically would conclude
Proposition~\ref{pr:bounded} from Proposition~\ref{pr:iss}. However,
the traditional ISS results require asymptotic convergence of the
unperturbed dynamics, (i.e., dynamics~\eqref{eq:ddyn} with $d \equiv
0$) to a point. This is not the case here and hence, 
we provide a formal proof.

\begin{remark}\longthmtitle{\baalgo is robust to variation in
    stepsizes}\label{re:robust-stepsizes}
  {\rm In practice, given that generators are competing and do not
    share information with each other, it is conceivable that they do
    not agree on a common stepsize. Propositions~\ref{pr:iss}
    and~\ref{pr:bounded} provide a way to quantify the performance of
    the algorithm when the stepsizes are different. Specifically, let
    $\beta_k$, $k \in \integerspositive$, denote a common set of
    stepsizes for all generators that satisfies the hypotheses of
    Theorem~\ref{th:convergence} and hence, guarantees the convergence
    properties outlined therein. Assume that each
    generator selects a different stepsize at each iteration, denoted as
    $\beta_{k,n}$, $k \in \integerspositive$, for generator $n$. Then,
    the bid iteration in Step~\ref{step3} of the \baalgo can be written
    as~\eqref{eq:ddyn} where now 
    \begin{align*}
      d_n(k) = (\beta_{k,n} - \beta_k)(x^{\opt}_n(k) - q_n(k))
    \end{align*}
    for all $n \in \until{N}$ and $k \in \integerspositive$. Now if the
    variation in stepsizes, i.e., the quantity $\beta_{k,n} - \beta_k$,
    is bounded above by a particular function of the distance of the
    bid-state to the efficient Nash equilibrium, then the linear 
    convergence and the
    ultimate bound is guaranteed following
    Proposition~\ref{pr:iss}. On the other hand, if the variation in
    stepsizes do not depend on the state but are bounded then, then the
    bids still converge asymptotically to a neighborhood of the
    efficient Nash equilibrium, as concluded in
    Proposition~\ref{pr:bounded}. Note that the assumption of $b_n(k)
    \ge c_n$ for all $n$ and $k$ holds whenever the stepsizes are
    positive for all agents at all times (cf.
    Lemma~\ref{le:seq-prop}).  \oprocend }
\end{remark}

\subsection{Robustness to deviation in bid update}

We illustrate here another aspect of robustness of the \baalgo by
establishing that, if all generators follow the bid update scheme,
then there is no incentive for any generator to deviate from it. We
next formalize these notions.  Assume that all generators, except $\nt
\in \until{N}$, follow the \baalgo, and that $\nt$ follows an
arbitrary strategy to update its bids. Then, one can write the \baalgo
under this deviation as
\begin{subequations}\label{eq:deviation-dyn}
  \begin{align}
    b_{-\nt}(k+1) & = [b_{-\nt}(k) + \beta_k (x^{\opt}_{-\nt}(k) - q_{-\nt}
    (k)) ]^+, \label{eq:deviation-dyn-1}
    \\
    b_{\nt} (k+1) & = \HH_{\nt}^{(k)} \Bigl( \{b_{\nt} (t), x^{\opt}_{\nt} (t),
    q_{\nt} (t)\}_{t=1}^k \Bigr),
    \label{eq:deviation-dyn-2}
    \\
    x^{\opt}(k+1) & \in \solone(b(k+1)), \label{eq:deviation-dyn-3}
    \\
    q(k+1) & \in \soltwo(b(k+1)), \label{eq:deviation-dyn-4}
  \end{align}
\end{subequations}
where the maps
$\{\map{\HH_{\nt}^{(k)}}{\realnonnegative^{3k}}{\realnonnegative}\}_{k=1}^\infty$
represent the update scheme of $\nt$ at iterations $1, 2,
\dots$ Recall that the subscript $-\nt$ denotes the vector without
the component corresponding to the generator $\nt$.  Note
that~\eqref{eq:deviation-dyn-2} implies that at each iteration $k$,
the generator $\nt$ only knows the bids it made and the quantities the
ISO demanded from it up until iteration $k$.

We next introduce the notion of ``incentive to deviate'' from the
\baalgo for the generator $\nt$.  A natural way to quantify incentives
for a generator is in terms of the payoff~\eqref{eq:payoff}: a
generator has an incentive to deviate if this would bring in a higher
payoff, when the ISO stops the iteration, than not deviating.  This is
formalized below.

\begin{definition}\longthmtitle{Incentive to deviate from
    \baalgo}\label{def:incentive-deviate}
  Let $r > 0$ and assume that the stepsizes for any execution
  of~\eqref{eq:deviation-dyn} satisfy the hypotheses of
  Theorem~\ref{th:convergence}. Then, the generator $\nt \in
  \until{N}$ has an \emph{incentive to deviate} from the \baalgo if
  there exists an execution of~\eqref{eq:deviation-dyn} and $l \in
  \integerspositive$ such that
    \begin{align}
      u_{\nt}(b_{\nt}(k), x_{\nt}^{\opt}(k)) > u_{\nt}^{\max},
      \label{eq:strong-dev-crit}
    \end{align}
    for all $k \ge l$, where 
    \begin{align} 
      u_{\nt}^{\max} := \max \Bigl\{u_{\nt}(b_{\nt}, &x_{\nt}^\opt(b))
      \Bigl|  \norm{b-b^*} \le \Bigl(1+\frac{B(r)}{2 a_{\max}} \Bigr) r
      \notag
      \\
      & \quad \text{ and } x^\opt (b) \in \solone(b)
    \Bigr\}. \label{eq:umax}
  \end{align}
\end{definition}

In the above definition, recall the short-hand notation $x^{\opt}(k)$
for $x^{\opt}(b(k))$. Equation~\eqref{eq:strong-dev-crit} implies that
the generator $\nt$ has an incentive to deviate if, after a finite
number of iterations, it is guaranteed a higher payoff than what it
might eventually get if it follows the \baalgo.  This captures the
fact that the generator does not know when the ISO might stop the bid
and hence it would deviate only when it is guaranteed to get a higher
payoff after a finite number of steps.
The next result shows that there is no incentive to deviate from the
\baalgo.

\begin{proposition}\longthmtitle{Robustness to
    deviation from \baalgo}\label{pr:robust-deviate-1}
  For dynamics~\eqref{eq:deviation-dyn}, let the hypotheses of
  Theorem~\ref{th:convergence} hold and assume that $b_n(k) \ge
  \coeffc_n$ for all $n \in \until{N}$ and $k \in \integerspositive$.
  Also, assume that the ISO selects a vertex solution $x^{\opt}(k) \in
  \solone(b(k))$ at each iteration $k \in \integerspositive$.  Then,
  no generator has an incentive to deviate from the \baalgo.
\end{proposition}
\begin{IEEEproof}
  We reason by contradiction. Assume that a generator $\nt$ has an
  incentive to deviate from the \baalgo. That is, there exists an
  execution of~\eqref{eq:deviation-dyn} and $l \in \integerspositive$
  such that~\eqref{eq:strong-dev-crit} holds for all $k \ge l$.  By
  definition,
  \begin{align}\label{eq:umax-lb-1}
    u^{\max}_{\nt} > b^*_{\nt} x^*_{\nt} - f_{\nt}(x^*_{\nt}).
  \end{align}
  Now consider the map 
  \begin{align*}
    \realnonnegative \ni b \mapsto g_{\nt}(b) := \max \setdef{ b q -
      f_{\nt}(q)}{q \ge 0}.
  \end{align*}
  From~\eqref{eq:efficiency}, we get $g_{\nt}(b^*_{\nt}) =
  b^*_{\nt} x^*_{\nt} - f_{\nt}(x^*_{\nt})$. Further,
  using~\eqref{eq:cost-quadratic}, one can show that this map is
  continuous, strictly increasing in the domain $b \ge c_{\nt}$, and
  $g_{\nt}(b) \to \infty$ as $b \to \infty$. These facts along
  with~\eqref{eq:umax-lb-1} imply that there exists a unique
  $b_{\nt}^{\max}
  > b^*_{\nt}$ such that $g_{\nt}(b_{\nt}^{\max}) = u_{\nt}^{\max}$,
  $g_{\nt}(b) > u_{\nt}^{\max}$ for all $b > b_{\nt}^{\max}$, and $g_{\nt}(b)
  < u_{\nt}^{\max}$ for all $\coeffc_{\nt} \le b < b_{\nt}^{\max}$.
  Then,~\eqref{eq:strong-dev-crit} reads as  
  \begin{align}\label{eq:umax-lb-2}
    u_{\nt}(b_{\nt}(k),x_{\nt}^{\opt}(k)) > g_{\nt}(b_{\nt}^{\max}),
  \end{align}
  for all $k \ge l$. From the above expression, we deduce that
  $b_{\nt}(k) \ge b_{\nt}^{\max}$ for all $k \ge l$.  Indeed
  otherwise, there exists $\tilde{k} \ge l$ such that
  $b_{\nt}(\tilde{k}) < b_{\nt}^{\max}$. This further implies that
  \begin{align*}
    u_{\nt}(b_{\nt}(\tilde{k}), x_{\nt}^{\opt}(\tilde{k})) 
    & = b_{\nt}(\tilde{k}) x_{\nt}^{\opt}(\tilde{k}) -
    f_{\nt}(x_{\nt}^{\opt}(\tilde{k}))
    \\
    & \le g_{\nt}(b_{\nt}(\tilde{k})) < g_{\nt}(b_{\nt}^{\max}),
  \end{align*}
  contradicting~\eqref{eq:umax-lb-2}. In the above expression, the first
  inequality follows from the definition of $g_{\nt}$ and the second
  follows from the fact that $g_{\nt}$ is strictly increasing. 

  The above reasoning has helped us establish that $b_{\nt}(k) \ge
  b_{\nt}^{\max} > b_{\nt}^*$ for all $k \ge l$. Note that
  $x_{\nt}^{\opt}(k) > 0$ for all $k \ge l$ because otherwise
  $u_{\nt}(b_{\nt}(k),x_{\nt}^{\opt}(k)) = 0$
  and~\eqref{eq:umax-lb-2} gets violated. By assumption, there exists
  at least one more generator connected to the bus $i(\nt)$ to which
  $\nt$ is connected to. For now assume that there is only one other
  generator $\nb \in \until{N}$ connected to $i(\nt)$. Since for all
  $k \ge l$, $x^{\opt}(k)$ is a solution of~\eqref{eq:sedp-flow},
  from the fact that $x_{\nt}^{\opt}(k) > 0$, we deduce
  \begin{align*}
    b_{\nb}(k) \ge b_{\nt}(k) \ge b_{\nt}^{\max},
  \end{align*}
  for all $k \ge l$. Now let
  \begin{align*}
    q_{\nb}^{\max} := \inf_{b \ge b_{\nt}^{\max}} \argmax \setdef{bq -
      f_{\nb}(q)}{q \ge 0}.
  \end{align*}
  Note that $q_{\nb}^{\max} > 0$ because of the facts: (i) $b^*_{\nb} =
  b^*_{\nt} < b^{\max}_{\nt}$; (ii) $\argmax \setdef{b^*_{\nb} q -
  f_{\nb}(q)}{q \ge 0} = x^*_{\nb} > 0$; and (iii) $b \mapsto \argmax
  \setdef{bq - f_{\nb}(q)}{q \ge 0}$ is nondecreasing.
  %
  %
  Since $b_{\nb}(k) \ge b_{\nt}^{\max}$ for all $k \ge l$, we obtain
  $q_{\nb}(k) \ge q_{\nb}^{\max}$ for all $k \ge l$ (see
  Step~\ref{step4} of the \baalgo for the definition of $q_{\nb}(k)$).
  Thus, if $b_{\nb} (k) > b_{\nt}(k)$ for some $k \ge l$, then
  $x_{\nb}^{\opt}(k) = 0$ (because $x^{\opt}(k)$ is an optimizer
  of~\eqref{eq:sedp-flow} given bids $b(k)$) and $b_{\nb}(k) >
  b_{\nt}^{\max}$.
  %
  %
  As a consequence,
  \begin{align}
    b_{\nb}(k+1) & = b_{\nb}(k) - \beta_k q_{\nb}(k)  
    \le b_{\nb}(k) - \alpha q_{\nb}^{\max}. \label{eq:finite-dec}
  \end{align}
  Therefore, if $b_{\nb}(k) > b_{\nt}(k)$ for some $k \ge l$, then
  from~\eqref{eq:finite-dec} we deduce that 
  there exists a finite $\tilde{k} > k$ such that, either
  $b_{\nb}(\tilde{k}) < b_{\nt}(\tilde{k})$ or $b_{\nb}(\tilde{k}) =
  b_{\nt}(\tilde{k})$. In the former case,
  $u_{\nt}(b_{\nt}(\tilde{k}),x_{\nt}^\opt (\tilde{k})) = 0$ as
  $x_{\nt}^{\opt}(\tilde{k}) = 0$. This
  contradicts~\eqref{eq:strong-dev-crit}.  In the latter case, two
  further cases can arise. In the first one, we get
  $x_{\nt}^{\opt}(\tilde{k}) = 0$ implying
  $u_{\nt}(b_{\nt}(\tilde{k}),x_{\nt}^{\opt}(\tilde{k})) = 0$ and
  contradicting~\eqref{eq:strong-dev-crit}. In the second one, we
  obtain $x_{\nb}^\opt (k) = 0$, implying $b_{\nb}(k+1) <
  b_{\nt}(k+1)$. This further yields
  $u_{\nt}(b_{\nt}(\tilde{k}+1),x_{\nt}^\opt (\tilde{k}+1)) = 0$,
  thereby, contradicting~\eqref{eq:strong-dev-crit}. Finally, if there
  are other generators connected to $i(\nt)$ that follow the \baalgo,
  then one can carry out the same reasoning as done above and show
  that we contradict~\eqref{eq:strong-dev-crit}. This completes the
  proof.
\end{IEEEproof}

\begin{remark}\longthmtitle{Generalization of
    Proposition~\ref{pr:robust-deviate-1}}\label{re:general-robust-deviate}
  {\rm It is interesting to observe that in the proof of
    Proposition~\ref{pr:robust-deviate-1}, we have not used at any
    point that the generators connected at buses other than the one
    that $\nt$ is connected follow the \baalgo. In fact, independently
    of how such generators update their bids, the \baalgo ensures that
    $\nt$ does not have any incentive to deviate. This is a useful
    property which we use later when studying robustness to collusion.
    \oprocend }
\end{remark}
\begin{remark}\longthmtitle{Other notions of ``incentive to
    deviate''}\label{re:other-incentives}
  {\rm In Definition~\ref{def:incentive-deviate}, one can impose the
    condition of higher payoff~\eqref{eq:strong-dev-crit} to hold for
    all executions of~\eqref{eq:deviation-dyn}. If this condition
    holds, then the generator has an even stronger incentive to
    deviate from the \baalgo. However, by
    Proposition~\ref{pr:robust-deviate-1}, we are ensured that there
    does not exist such strong incentive to deviate. This is because
    the result shows that there does not exist any execution
    of~\eqref{eq:deviation-dyn} for which~\eqref{eq:strong-dev-crit}
    holds. If, on the other hand, we replace the
    condition~\eqref{eq:strong-dev-crit} in
    Definition~\ref{def:incentive-deviate} with the requirement that
    there exists an execution of~\eqref{eq:deviation-dyn} along which
    \begin{align}
      \limsup_{k \to \infty}  u_{\nt} (b_{\nt}(k), x_{\nt}^{\opt}(k))
      > u_{\nt}^{\max} \label{eq:weak-dev-crit}
    \end{align}
    holds.  This inequality means that there exists an execution
    of~\eqref{eq:deviation-dyn} in which the generator $\nt$ gets a
    higher payoff than $u_{\nt}^{\max}$ infinitely often. Since the
    ISO can stop the iterations at any time, the generator is not
    guaranteed a higher payoff, but the possibility is still there. We
    conjecture that the \baalgo is not robust to this notion of weak
    incentive to deviate. However, the obfuscation of the stopping
    criteria by the ISO makes such a weak incentive not enough for a
    rational generator to deviate.  \oprocend }
\end{remark}

\subsection{Robustness to collusion}\label{sec:collusion}

Here we study the robustness of the \baalgo against
collusion. Collusion refers to the action of a set of generators
to share among themselves information about their bids and generation
demands by the ISO, with the goal of getting a higher profit,
possibly by deviating from the bid update scheme.
The following makes this notion formal.

\begin{definition}\longthmtitle{Collusion between
    generators}\label{def:collusion}
  A group of generators $\JJ \subset \until{N}$ form a collusion if at
  each iteration $k \in \integerspositive$ of the algorithm, each
  generator $n \in \JJ$,
  \begin{enumerate}
  \item has the information 
    \begin{align*}
      \II_k := \setdef{(b_{r}(t),x^{\opt}_{r}(t))}{r \in \JJ, t
        \in \until{k}}, \text{ and }
    \end{align*}
  \item determines its next bid $b_{n}(k+1)$ based on the information
    $\II_k$, not necessarily following the update scheme
    (Step~\ref{step3}) of the \baalgo.
    \end{enumerate}
\end{definition}


An iteration of the \baalgo under a collusion between a group of
generators $\JJ \subset \until{N}$ is
given by the following dynamics
\begin{subequations}\label{eq:collusion-dyn}
    \begin{align}
      b_n(k+1) & = [b_n(k) + \beta_k (x_n^{\opt}(k) - q_n(k))]^+, \forall n 
      \not \in \JJ,
      \label{eq:collusion-dyn-1} 
       \\
       b_{n}(k+1) & =
       \HH_{n}^{(k)}\Bigl(\II_k, \{q_{n}(t)\}_{t=1}^k \Bigr), \forall n \in \JJ
       \label{eq:collusion-dyn-2}
       \\ 
       x^{\opt}(k+1) & \in \solone(b(k+1)), \label{eq:collusion-dyn-3}
       \\
       q(k+1) & = \soltwo(b(k+1)), \label{eq:collusion-dyn-4}
    \end{align}
\end{subequations}
where maps $\{\map{\HH_{n}^{(k)}}{\realnonnegative^{(2
    \abs{\JJ}+1)k}}{\realnonnegative}\}_{n \in \JJ, k=1, 2, \dots}$
represent the update scheme of generators in collusion.  Notice that
for each generator $n$, the quantity $q_n(k)$, for all $k \in
\integerspositive$, is part of its private information, irrespective
of the fact that $n$ belongs to $\JJ$ or not.  Next, we define what it
means for the group of generators $\JJ$ to have an incentive to
collude.

\begin{definition}\longthmtitle{Incentive to
    collude}\label{def:incentive-collude}
  Let $r > 0$ and assume that the stepsizes for any execution
  of~\eqref{eq:collusion-dyn} satisfy the hypotheses of
  Theorem~\ref{th:convergence}. Then, the group of generators $\JJ$
  has an \emph{incentive to collude} under the \baalgo if there exists
  an execution of~\eqref{eq:collusion-dyn}, a generator $\nt \in \JJ$,
  and $l \in \integerspositive$ such that
    \begin{align}
      u_{\nt}(b_{\nt}(k), x_{\nt}^{\opt}(k)) > u_{\nt}^{\max},
      \label{eq:collude-dev-crit}
    \end{align}
    for all $k \ge l$, where $u_{\nt}^{\max}$ is
    defined in~\eqref{eq:umax}.
\end{definition}

This notion essentially says that there is an incentive to collude for
the generators in $\JJ$ if there exists at least one execution
of~\eqref{eq:collusion-dyn} along which at least one generator in
$\JJ$ gets a higher payoff after finite number of
steps. The next result shows that no group of generators has an
incentive to collude provided there is at least one generator at each
bus with generation that follows the \baalgo.


\begin{proposition}\longthmtitle{Robustness to
    collusion under the \baalgo}\label{pr:robust-collude}
  For dynamics~\eqref{eq:collusion-dyn}, let the hypotheses of
  Theorem~\ref{th:convergence} hold and assume that $b_n(k) \ge
  \coeffc_n$ for all $n \in \until{N}$ and $k \in \integerspositive$.
  Assume that the ISO selects a vertex solution $x^{\opt}(k) \in
  \solone(b(k))$ at each iteration $k \in \integerspositive$.  Assume
  that at each bus that has generators connected to it, there exists
  at least one generator that follows the update scheme of the
  \baalgo. Denote these generators by $\KK \subset \until{N}$. Then,
  there is no incentive to collude for any group of generators
  contained in $\until{N} \setminus \KK$.
\end{proposition}
\begin{IEEEproof}
  Let $\JJ \subset \until{N} \setminus \KK$ be a group of generators
  that form a collusion. Assume first Scenario~$1$ 
  where each generator in $\JJ$ is connected to a different bus.
  By hypotheses, there exists at least one other generator following
  the \baalgo at the bus where a generator in $\JJ$ is connected
  to. Thus, mimicking the proof of
  Proposition~\ref{pr:robust-deviate-1} (cf.
  Remark~\ref{re:general-robust-deviate}), at each bus, no generator
  has an incentive to deviate from the $\baalgo$.  By
  Definition~\ref{def:incentive-deviate}, this implies that there does
  not exist any execution of~\eqref{eq:collusion-dyn} for
  which~\eqref{eq:collude-dev-crit} holds for any generator in~$\JJ$.
  Hence, for Scenario $1$, generators in $\JJ$ do not have an
  incentive to collude.

  Next, consider Scenario $2$, where at least a bus, say $i \in
  \until{N_b}$, has more than one generator from $\JJ$, that is,
  $\JJ_i : = G_i \cap \JJ$ has cardinality larger than or equal to
  $2$. Let $\nb \in G_i$ be the generator at $i$ that follows
  \baalgo. For the sake of contradiction, assume the existence of a
  generator $\nt \in \JJ_i$ for which~\eqref{eq:collude-dev-crit}
  holds for some execution of~\eqref{eq:collusion-dyn}. Since the ISO
  selects a vertex solution at each iteration $k \in
  \integerspositive$, we deduce that for all $k \ge l$, all other
  generators in $\JJ_i$ get zero production signal from the ISO, i.e.,
  $x^{\opt}_n(k) = 0$ for all $n \in \JJ_i \setminus \{\nt\}$ and $k
  \ge l$. Therefore, for the purpose of analysis, one can neglect the
  generators in $\JJ_i \setminus \{\nt\}$ and assume that only $\nt$
  and $\nb$ are connected to $i$. Again, mimicking the proof of
  Proposition~\ref{pr:robust-deviate-1}, we deduce that $\nt$ does not
  have an incentive to deviate and so~\eqref{eq:collude-dev-crit} does
  not hold, a contradiction. Since $i$ is arbitrary, we conclude that
  for Scenario $2$, generators in $\JJ$ do not have an incentive to
  collude either.
  %
  %
\end{IEEEproof}

An alternative definition of an incentive to collude
could be where every generator in the collusion gets a higher
payoff after a finite number of steps.
Proposition~\ref{pr:robust-collude} however shows that,
under the assumed hypotheses, such a scenario does not occur as there
is not even a single generator that gets a higher payoff after a
finite number of iterations. 
Note that the assumptions of the above result is tight in the sense that
if all generators at a bus collude, then based on the load and the line
limits, generators at that bus can increase their bid
to an arbitrarily high value, thus creating an incentive to collude.

\begin{remark}\longthmtitle{Limitations on robustness under generator
  bounds}\label{re:limitations}
{\rm The robustness of the \baalgo against deviation and collusion
  relies heavily on the fact that we have not considered upper bounds
  on the generation capacities. In the presence of such bounds, the
  generators might be able to push the bids and their individual
  utilities to a higher value based on the load at the respective bus
  and the capacity constraints on the lines connected to the bus. To
  avoid such behavior of market manipulation, either one can modify
  network capacities or investigate alternative allocation mechanisms
  that disincentivizes such behavior.  \oprocend }
\end{remark}

\section{Simulations}\label{sec:simulations}
\begin{figure}[htb]
  \centering
  \includegraphics[width=.68\linewidth]{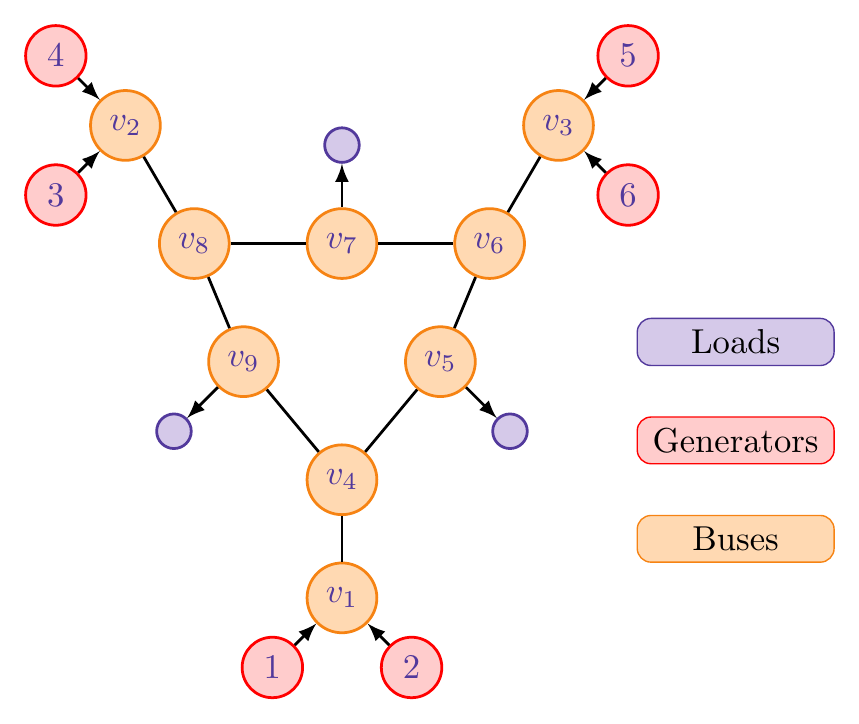}
  \caption{Network layout of the modified IEEE $9$-bus test
    case.}\label{fig:network}
    \vspace*{-2ex}
\end{figure}

\begin{figure}[htb]
  \centering
  \subfigure[Bids
  $b(k)$]{\includegraphics[width=0.42\linewidth]{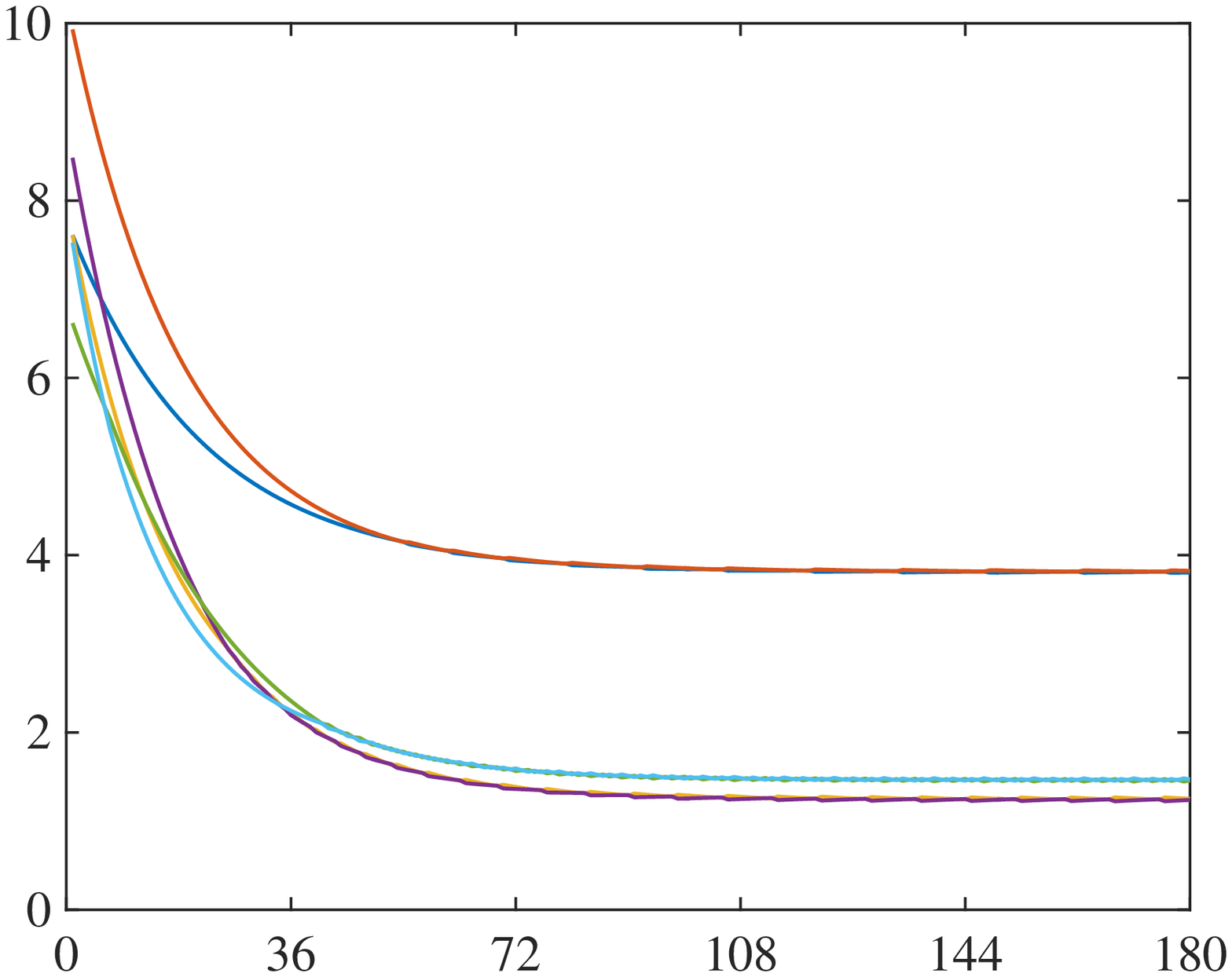}}
  \, \, 
  \subfigure[Distance
  $\norm{b(k)-b^*}$]{\includegraphics[width=0.42\linewidth]{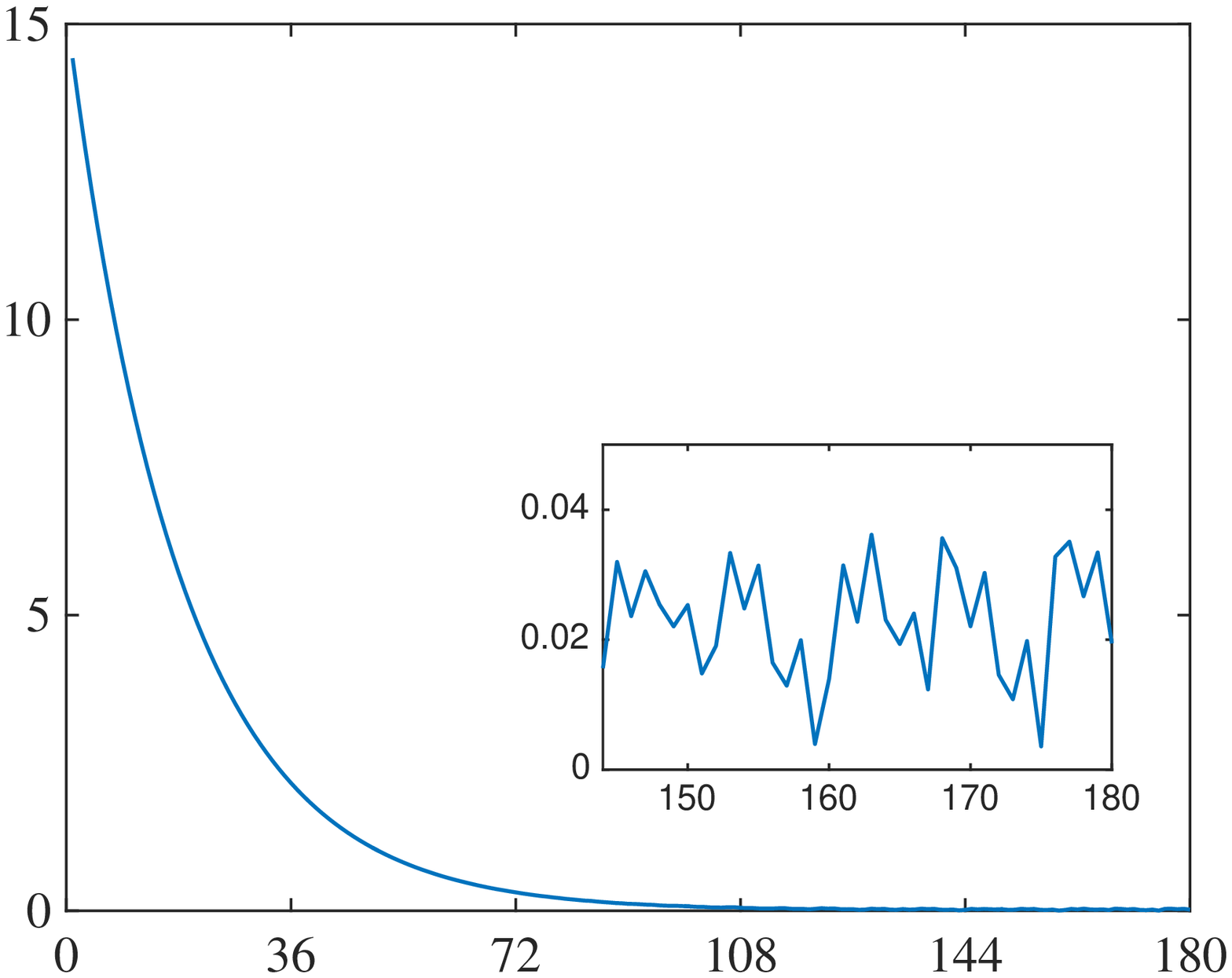}} 
  \, \,
  \caption{Execution of the \baalgo for the modified IEEE 9-bus test
    case in Figure~\ref{fig:network}. The cost function for each
    generator $i$ is $f_i(x_i) = a_i x_i^2 + c_i x_i$, with
    coefficients given in~\eqref{eq:coeffs}. The load is $y_5 = 2$,
    $y_7 = 3$, and $y_9 = 1$. The efficient Nash equilibrium $b^*$ is
    given in~\eqref{eq:effNE}. Plots (a) and (b) show, respectively,
    the evolution of the bids and their distance to~$b^*$. The
    stepsizes are $\beta_k = 0.01$ for all $k$ and the initial bids
    are $b(1) = (7.6096, 9.9313, 7.6087, 8.4827, 6.6175, 7.5254)$.
    Bids converge to and then remain in a neighborhood of the
    efficient Nash equilibrium.}
    \label{fig:evol}
    \vspace*{-2ex}
\end{figure}

We illustrate the convergence and robustness properties of the \baalgo
using a modified IEEE $9$-bus test case~\cite{RDQ-CEMS-RJT:11}.  The
traditional IEEE 9-bus system has 3 generators, at buses $v_1$, $v_2$,
and $v_3$ and three loads at buses $v_5$, $v_7$, and $v_9$. 
In our modified test case, we have added one generator each at buses
$v_1$, $v_2$ and $v_3$. The interconnection topology is given in
Figure~\ref{fig:network}. 
%
%
The line flow limit between any two buses $(v_i,v_j)$ is $2.5$ except
for three lines, $(v_5,v_6)$, $(v_3,v_6)$, and $(v_6,v_7)$, for which
the limits are $1.5$, $3.0$, and $1.5$, respectively.  The loads are
$y_5 = 2$, $y_7 = 3$, and $y_9 = 1$, where $y_i$ denotes the load at
bus $v_i$.  The cost function for each generator $i$ is $f_i(x_i) =
a_i x_i^2 + c_i x_i$, where the coefficients for all the generators are
given by the vectors
\begin{align}
  a & = (0.1100, 0.0950, 0.0850, 0.1000, 0.1225, 
   0.0750),\notag 
  \\
  c & = (3.5, 3.8, 1.2, 0.8, 1.0, 
  1.3).  \label{eq:coeffs} 
\end{align}
For the given costs and loads, the generation profile at the optimizer
of the DC-OPF problem~\eqref{eq:edp-flow} is 
\begin{align*}
  x^* & = (1.4268, 0.0732, 0.2703, 2.2297, 1.8987, 1.1013),
\end{align*}
and the unique efficient Nash equilibrium is 
\begin{align}
  b^* & = (3.8139, 3.8139, 1.2459, 1.2459, 1.4652, 1.4652).   \label{eq:effNE}
\end{align}
Figure~\ref{fig:evol} depicts the evolution of the bids and their
distance to the efficient Nash equilibrium along an execution of the
\baalgo. The initial bids $b(1)$ are selected satisfying 
$b_n(1) \ge c_n$ for all the generators $n \in \until{6}$.
The stepsizes are constant, $\beta_k = 0.01$ for all $k$, and satisfy
 $\beta_k < 2 a_n$.  As predicted by
Theorem~\ref{th:convergence}, Figure~\ref{fig:evol} shows that the
bids converge towards the efficient Nash equilibrium $b^*$ at a linear
rate and, after a finite number of steps, remain in a neighborhood
of~$b^*$.  If one selects $r = 1.35$, then $B(r) = 0.0101$ and
condition~\eqref{eq:beta-cond-3} holds for the stepsizes.
Computing the right hand side of~\eqref{eq:ultimate-bound} using these
values, we conclude that bids eventually remain in the neighborhood
centered at $b^*$ with radius $1.3775$.  Figure~\ref{fig:evol}(b)
validates this claim and, in fact, shows that the bound is
conservative since bids actually remain in a neighborhood of
radius~$0.05$.
 
\begin{figure}[htb]
  \centering
  \subfigure[Bids
  $b(k)$]{\includegraphics[width=0.42\linewidth]{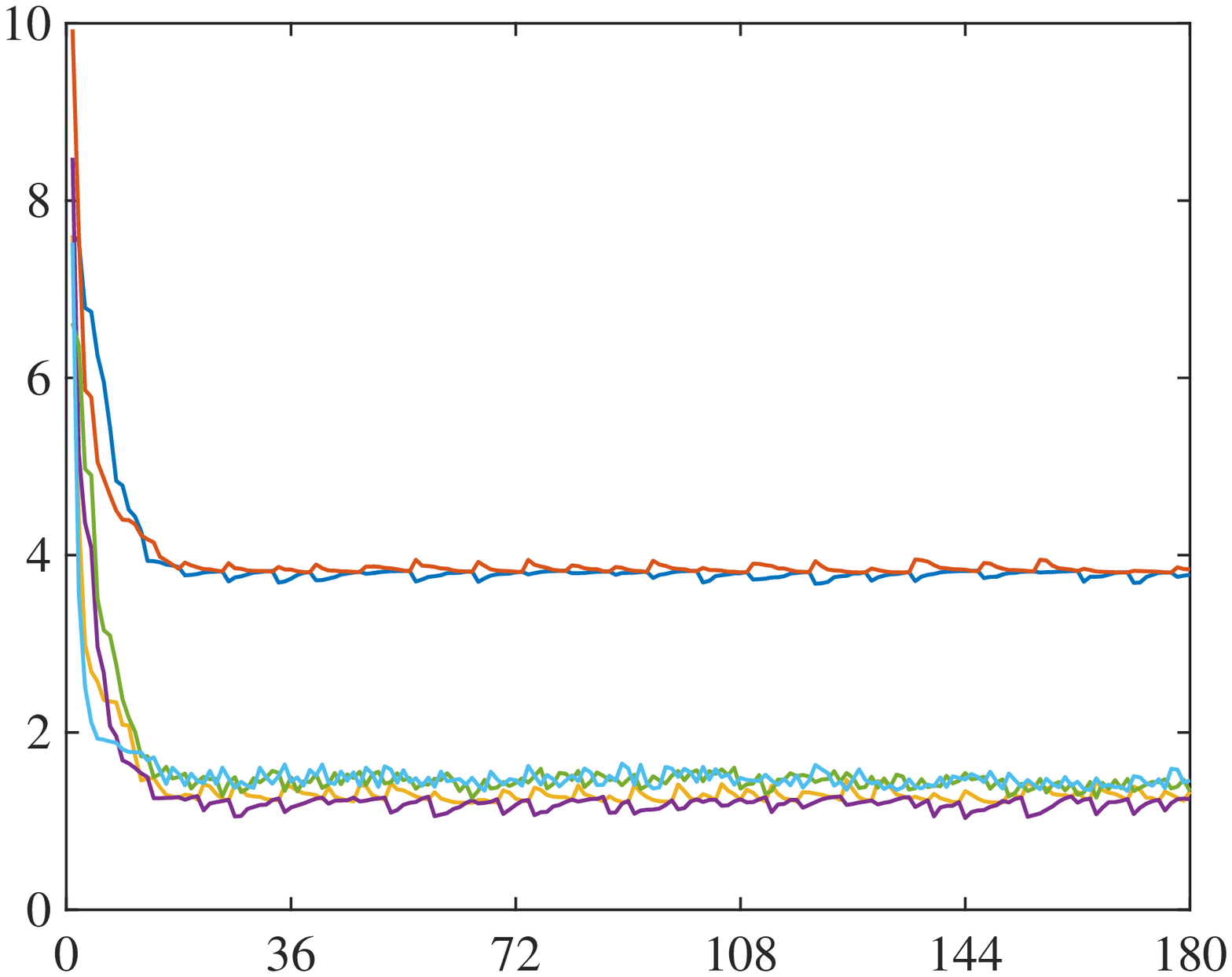}}
  \, \, 
  \subfigure[Distance
  $\norm{b(k)-b^*}$]{\includegraphics[width=0.42\linewidth]{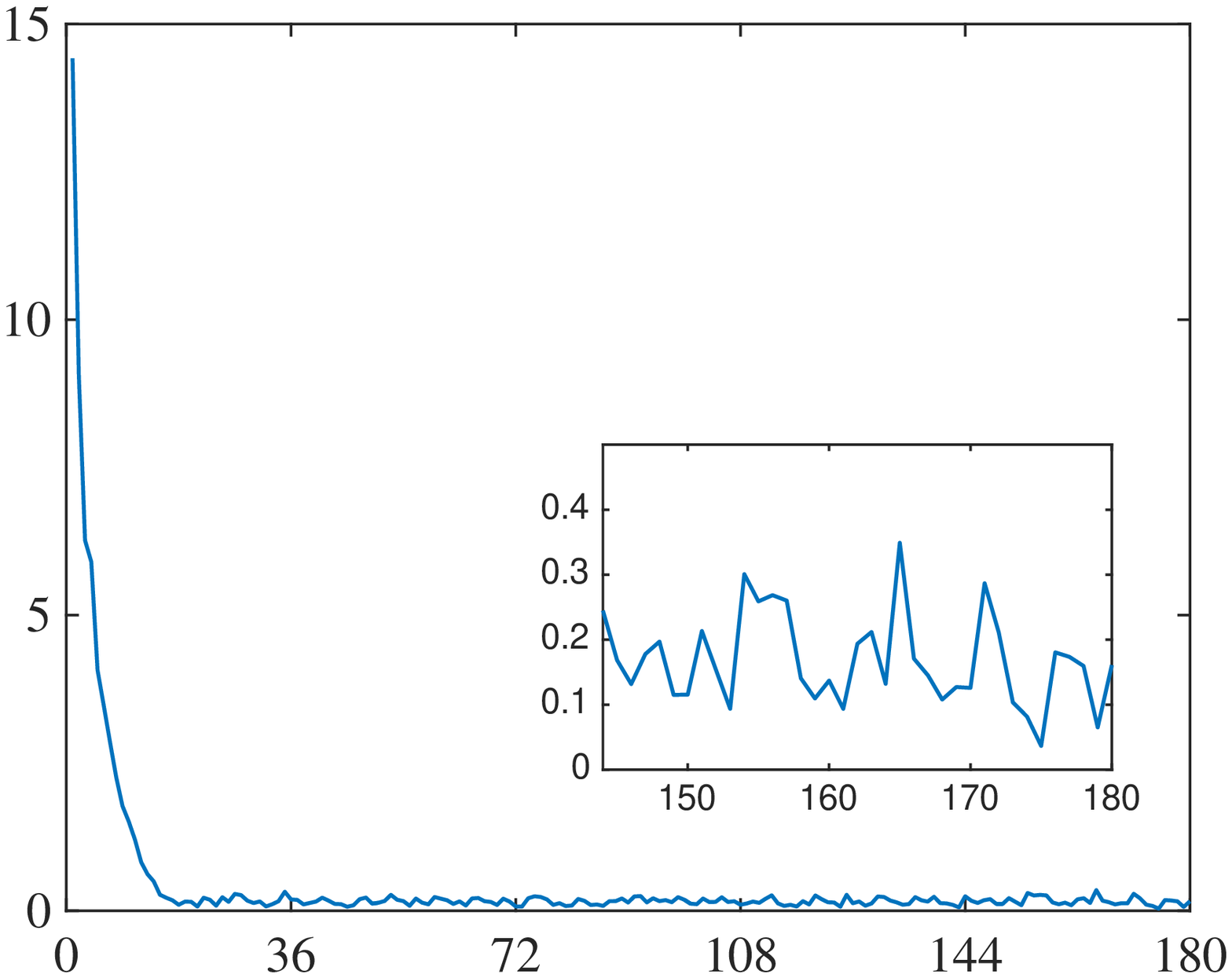}} 
  \, \,
  \\
  \subfigure[Bids
  $b(k)$]{\includegraphics[width=0.42\linewidth]{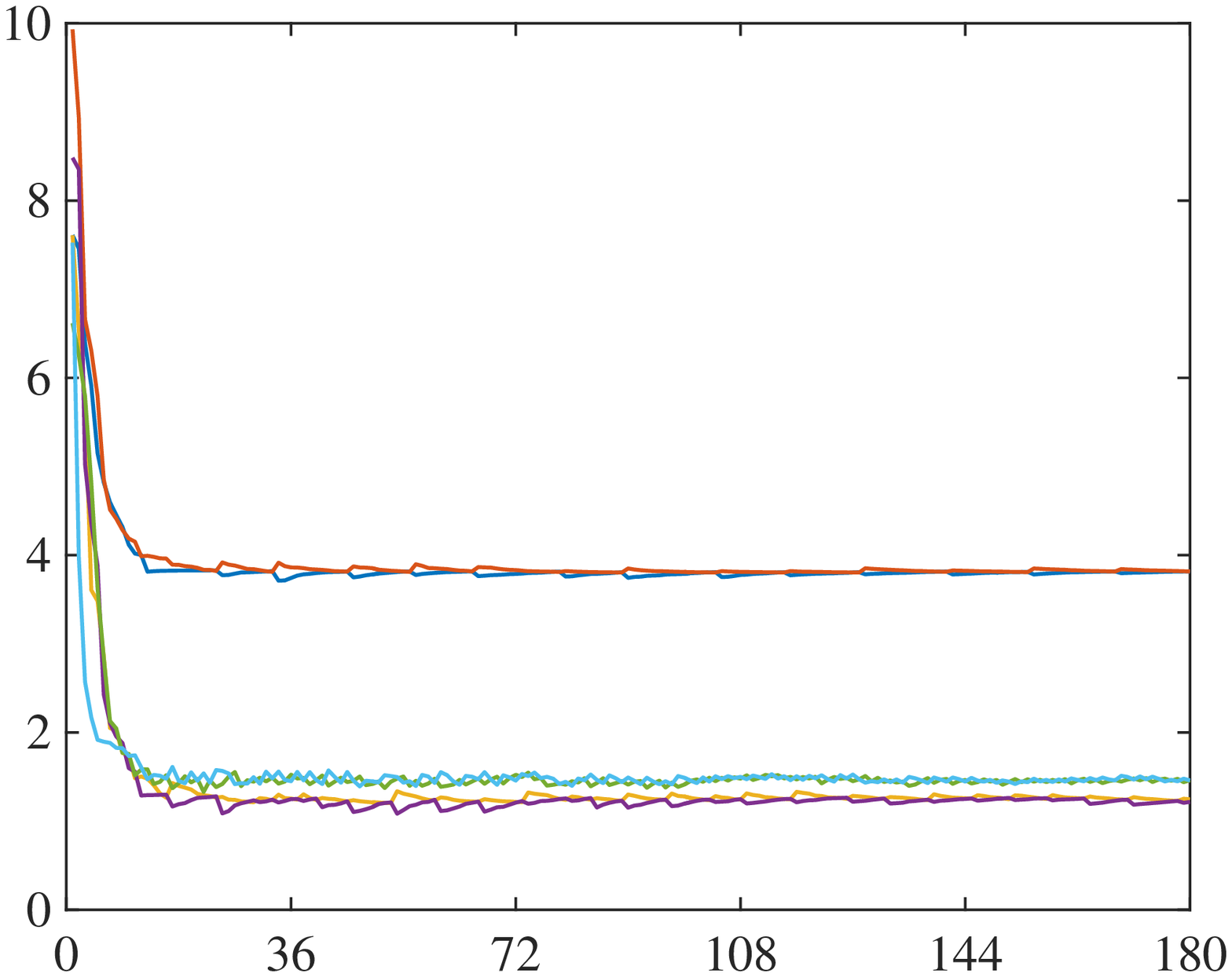}}
  \, \, 
  \subfigure[Distance
  $\norm{b(k)-b^*}$]{\includegraphics[width=0.42\linewidth]{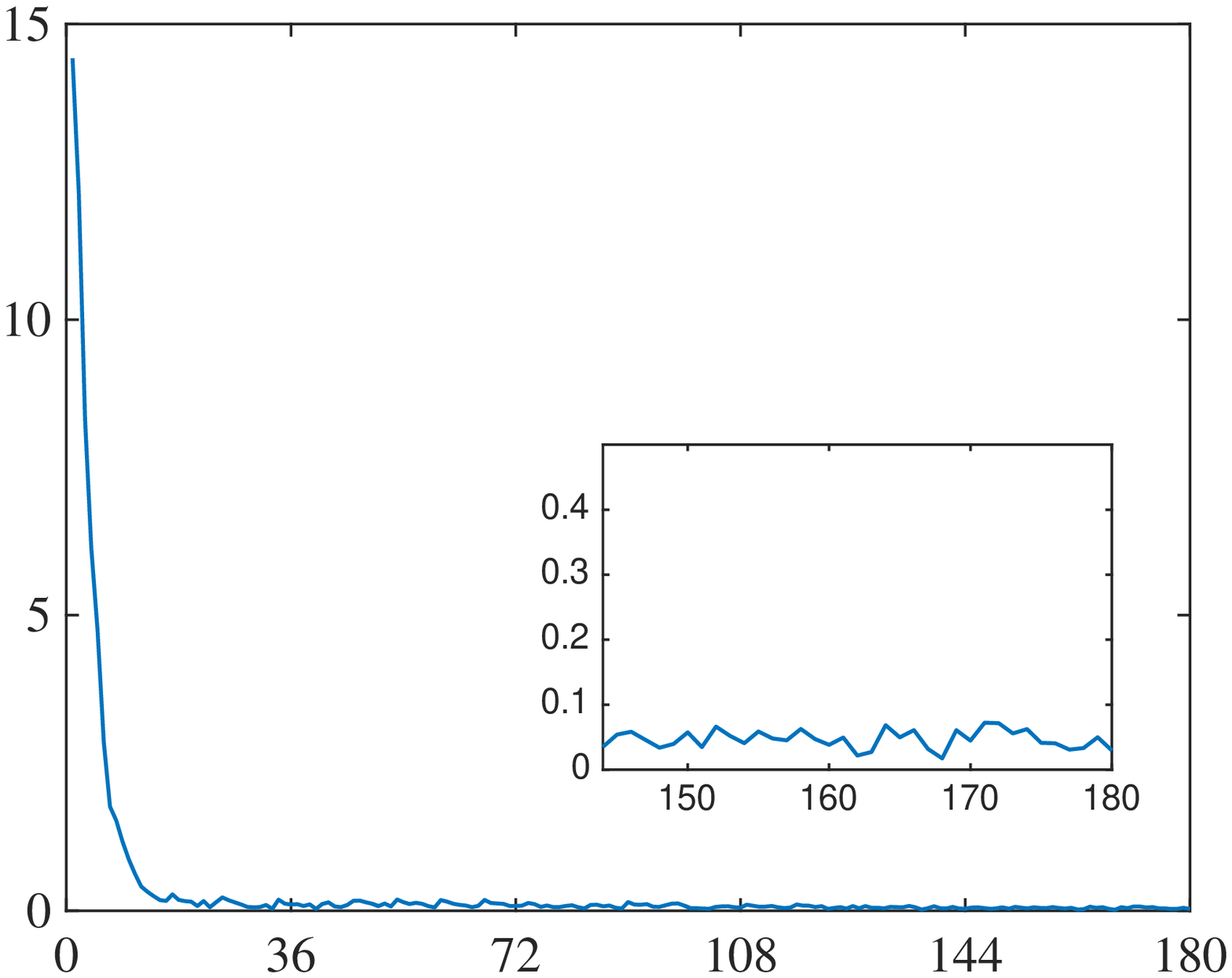}} 
  \, \,
  \caption{Execution of the \baalgo under different stepsize selection
    for the example of Figure~\ref{fig:evol}. All data is the same
    except for the stepsizes. In plots (a) and (b), each generator at
    each iteration randomly selects the stepsize from the set
    $[0.001,0.1]$ with uniform probability distribution. We observe
    that the bids still converge to a neighborhood of the efficient
    Nash equilibrium, but the size of the neighborhood is bigger than
    that achieved in Figure~\ref{fig:evol}. In plots (c) and (d), the
    interval of stepsize selection decays with time to a single point
    $0.01$. The bids now converge to the efficient Nash equilibrium
    with greater accuracy. These observations validate the robustness
    guarantees of Proposition~\ref{pr:bounded}.}
    \label{fig:dist}
    \vspace*{-2ex}
\end{figure}

\begin{figure}[htb]
  \centering
  \subfigure[$u_n(b_n(k),x^{\opt}_n(k))-u_n(b^*_n,x^*_n)$]{\includegraphics[width
  = 0.55 \linewidth]{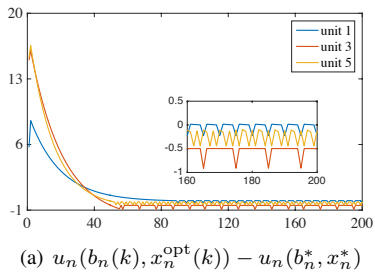}}
\caption{Execution of the \baalgo for the example considered in
  Figure~\ref{fig:evol} with generators $1$, $3$, and $5$ forming a
  collusion. The initial condition is the same and the stepsize is
  $0.01$ at each iteration for generators $2$, $4$, and $6$.  For each
  $n \in \{1,3,5\}$, at each iteration $k$, $b_n(k) = 0.99*b_{n+1}(k)$
  if this value is bigger than or equal to $b^*_n$. Otherwise,
  $b_n(k)$ is selected randomly from the interval $[b_n^*,b_n^*+1]$,
  with uniform probability distribution.  With this choice of bid, the
  colluding generators aim to get a positive production signal and at
  the same time bid high enough so as to obtain a high utility.  The
  plot shows the evolution of the difference between the utility
  obtained at each iteration, $u_n(b_n(k),x^{\opt}_n(k))$, and the
  utility at the optimal bid and generation,
  $u_n(b_n^*,x^{\opt}_n(b^*))$ for each $n \in \{1,3,5\}$.  This value
  becomes negative for all generators after a finite number of
  iterations. Since $u_n^{\max} > u_n(b_n^*,x^*_n)$, the example shows
  that~\eqref{eq:collude-dev-crit} does not
  hold.}\label{fig:collusion}
  \vspace*{-2ex}
\end{figure}

We next illustrate the robustness properties of the \baalgo against
disturbances (cf. Section~\ref{sec:robust-dist}).
Figure~\ref{fig:dist} considers the same setup as above but now with
generators choosing a different stepsize at each iteration. These
differences in stepsizes can be interpreted as a disturbance to the
\baalgo, as discussed in Remark~\ref{re:robust-stepsizes}. In
Figure~\ref{fig:dist}(a)-(b), the interval from which stepsizes are
selected is constant, whereas in Figure~\ref{fig:dist}(c)-(d) the size
of this interval decays with time. In both cases, the bids converge to
a neighborhood of $b^*$ (in the latter case
of decaying interval, the bids converge to a smaller neighborhood), as
established in Proposition~\ref{pr:bounded}. Observe that the
convergence rate in Figure~\ref{fig:dist}(a)-(b) is higher than in
Figure~\ref{fig:evol}(a)-(b). This is because 
stepsizes are allowed to be large in the former. However, 
this higher convergence rate comes with the pitfall of loss in
accuracy, cf. Remark~\ref{re:conv-prop}. Hence, to retain both
properties, stepsizes should be large initially and decay as
iterations proceed. This is seen in Figure~\ref{fig:dist}(c)-(d),
where stepsizes decay over time (in expectation), yielding
both high convergence rate and accuracy.  Finally,
Figure~\ref{fig:collusion} demonstrates the robustness against
collusion of the \baalgo (cf.  Section~\ref{sec:collusion}), where
generators $1$, $3$, and $5$ form a collusion.  These generators may
select their bids in any fashion they want: for this example, we
assume a particular strategy of bid selection, explained in
Figure~\ref{fig:collusion}.
The plot shows that the utility of the colluding generators eventually
becomes lower than $u_n^{\max}$ (defined in~\eqref{eq:umax}). Hence,
there is no incentive for collusion, as ensured by
Proposition~\ref{pr:robust-collude}.


\section{Conclusions}\label{sec:conclusions}

We have formulated an inelastic electricity market game capturing the
strategic interaction between generators in a bid-based energy
dispatch setting. For this game, we have established the existence and
uniqueness of the efficient Nash equilibria. We have also designed the
\baalgo, which is an iterative strategy amenable to decentralized
implementation that provably converges to a neighborhood of the
efficient Nash equilibrium at a linear rate.  We have characterized
the robustness properties of the algorithm against disturbances,
deviation in bid updates, and collusions among generators. Future work
will analyze the dynamic behavior of the market under other bidding
schemes, such as Cournot bidding, supply function bidding, and
price-capacity bidding.  We would also like to examine the convergence
of other learning schemes such as regret minimization in the context
of electricity markets.  Finally, we wish to incorporate stochastic
load demands and changing sets of generators in our setup.
 

\bibliographystyle{IEEEtran}%
\bibliography{alias,Main,Main-add,JC}

\end{document}